\documentclass{amsart}

\usepackage{amsmath}
\usepackage{amsthm,amssymb,color,comment,yhmath}
\usepackage{amsthm,amssymb,color,comment}
\usepackage{mathtools}
\usepackage{bbm}
\usepackage{hyperref}
\usepackage[TS1,T1]{fontenc}
\usepackage{inputenc}
\usepackage{dsfont}
\usepackage{tikz}
\usepackage{enumitem}
\usepackage{soul}
\usepackage{setspace}
\usepackage[top = 2.5cm, left=2cm, right=2cm, bottom=2.5cm]{geometry}
\usepackage[sort]{cite}

\numberwithin{equation}{section}

\newtheorem{thm}{Theorem}[section]
\newtheorem{prop}[thm]{Proposition}
\newtheorem{lem}[thm]{Lemma}

\newtheorem{Def}[thm]{Definition}

\theoremstyle{definition}

\newtheorem{rem}[thm]{Remark}

\DeclareMathOperator{\DIV}{div}

\newcommand{\R}{\mathbb{R}}

\newcommand{\T}{\mathbb{T}}
\newcommand{\p}{\partial}
\newcommand{\I}{\mathbb{Z}}
\newcommand{\eps}{\varepsilon}

\newcommand{\diff}{\mathop{}\!\mathrm{d}}

\makeatletter
\newcommand{\doublewidetilde}[1]{{%
  \mathpalette\double@widetilde{#1}%
}}
\newcommand{\double@widetilde}[2]{%
  \sbox\z@{$\m@th#1\widetilde{#2}$}%
  \ht\z@=.9\ht\z@
  \widetilde{\box\z@}%
}
\makeatother

\author{Piotr Gwiazda}
\address{{\it Piotr Gwiazda:} Institute of Mathematics, Polish Academy of Sciences, ul. Śniadeckich 8, 00-656 Warsaw, Poland; Interdisciplinary Centre for Mathematical and Computational Modelling, University of Warsaw, ul. Tyniecka 15/17, 02-630 Warsaw, Poland}
\email{pgwiazda@mimuw.edu.pl}

\author{Jakub Skrzeczkowski}
\address{{\it Jakub Skrzeczkowski: } {St John's College, University of Oxford, St Giles, Oxford, OX1 3JP, United Kingdom} \& Mathematical Institute, University of Oxford, Woodstock Road, Oxford, OX2 6GG, United Kingdom}
\email{jakub.skrzeczkowski@maths.ox.ac.uk}

\author{Jakub Woźnicki}
\address{Institute of Mathematics of Polish Academy of Sciences, Jana i J\k edrzeja \'Sniadeckich 8, 00-656 Warsaw, Poland; Faculty of Mathematics, University of Vienna}
\email{jakub.woznicki@univie.ac.at}

\begin{document}

\title[Convergence rates for smoothed particle hydrodynamics (SPH) with nonlinear pressure]{Convergence rates for smoothed particle hydrodynamics (SPH)\\ with nonlinear pressure}

\begin{abstract}
We establish quantitative convergence results for the smoothed particle hydrodynamics (SPH) method approximating the compressible Euler equations with nonlinear pressure laws $\rho^\gamma$, $\gamma\geq 2$. Our approach is based on a nonlocal approximation of the pressure force that preserves the energy structure of the original system, including an exact energy balance. Using the relative entropy method, we prove convergence towards smooth solutions of the compressible Euler equations as the nonlocality parameter $\eps\to0$, with algebraic convergence rates and an explicit algebraic relationship between $\eps$ and the number of particles $N$. In particular, we obtain quantitative convergence rates for the density in the Wasserstein distance and negative Sobolev norms, and for the momentum in the flat metric. A key feature of our analysis is that it applies to measure-valued densities in the nonlocal system, providing a unified framework for both continuum and particle solutions. Our results extend the existing rigorous convergence theory for SPH beyond quadratic pressure ($\gamma=2$) to general nonlinear pressure laws with $\gamma\geq 2$ on the periodic domain.
\end{abstract}

\keywords{smoothed particle hydrodynamics (SPH), compressible Euler equations, nonlinear pressure, nonlocal-to-local convergence, relative entropy method, particle approximations, quantitative convergence rates, Wasserstein distance}
\subjclass[2020]{65M15, 65M12, 35Q31, 76M28, 65M75, 35B35}
\maketitle

\section{Introduction}

\noindent In this paper we study the smoothed particle hydrodynamics method (SPH) for the compressible Euler equation
\begin{equation}\label{sys:euler_general}
\begin{aligned}
    \partial_t\rho + \DIV_x(\rho\, u) &= 0,\\
    \partial_t (\rho\,u)  + \DIV_x (\rho\, u \otimes u ) + \rho\nabla_x \rho^{\gamma-1}+ \kappa \,\rho\, u  &= 0. \\
    \end{aligned}
\end{equation}
This method, introduced in \cite{gingold1977smoothed, lucy1977numerical, MR673709} in the context of numerical simulations of self-gravitating fluids in astrophysics, is based on approximating the nonlinear pressure term via convolution. In our case, we will be approximating 
$$
\rho^{\gamma-1} \approx (\rho \ast L_{\eps})^{\gamma-1} \ast L_{\eps},
$$
 where $L_{\eps}$ is a symmetric smoothing kernel approximating the Dirac delta as $\eps\to0$. This particular form respects the energy structure of the equation and is crucial for our approach to succeed but other forms are possible. The nonlocal approximation allows us to solve the equation using a Lagrangian approach. Assuming the fluid is composed of $N$ particles with positions $X_1(t), \ldots, X_N(t)$ and velocities $V_1(t), \ldots, V_N(t)$, they satisfy the ODEs:
\begin{equation*}
    \begin{split}
        \dot{X_i}(t) &= V_i(t),\\
        \dot{V_i}(t) &= -\int_{\R^d}\nabla_x L_{\eps}(X_i(t) - y)\left(\sum_{k=1}^N m_k\, L_{\eps}(y - X_k(t)) \right)^{\gamma-1}\diff y - \kappa \,V_i(t). 
    \end{split}
\end{equation*}
The SPH method remains widely used in computational physics, with recent applications to astrophysical simulations of galaxy formation and turbulence, see \cite{schaller2024swift} and \cite{cabezon2026particles}, respectively. For broader overviews of its applications to astrophysical systems, including compact objects, and to free-surface and multiphase flows, we refer to the reviews \cite{rosswog2015sph} and \cite{le2025smoothed}, respectively.\\

\noindent The properties of the SPH method are closely related to those of the corresponding Euler system, i.e. when $\rho^{\gamma-1}$ is approximated via convolution. Indeed, the idea is to separate the convergence analysis into two steps: first, establishing convergence of the particle system to the nonlocal Euler equations as $N\to\infty$ for fixed $\varepsilon$, and then studying the limit $\varepsilon\to0$ in the nonlocal equations. This approach was followed by Di Lisio, Grenier, and Pulvirenti \cite{MR1605143}, who established convergence of SPH approximations to a regularized Euler system using a Lagrangian formulation and estimates in the Wasserstein distance (see also \cite{MR3756760} for a similar reasoning). In their approach, however, the stability estimates involve constants that grow exponentially in negative powers of the nonlocality parameter $\eps$ (these are norms of derivatives of nonlocal operator which naturally blow up when $\eps \to 0$), leading to a logarithmic relationship between $\varepsilon$ and the number of particles $N$ when both limits are taken simultaneously, which is rather restrictive in practice. Concerning the nonlocal equations themselves, Di Lisio, Grenier, and Pulvirenti \cite{MR1487955} established uniform estimates in Sobolev spaces for nonlinear pressure laws on a sufficiently short time interval. More recently, Danchin and Mucha \cite{MR4764519} obtained uniform estimates in Besov spaces, although their analysis is restricted to quadratic pressure. In both cases, these estimates imply convergence of sequence of solutions $\rho_{\eps}$ of the nonlocal equation to a solution of the local Euler system \eqref{sys:euler_general} as $\varepsilon\to0$. The problem of the nonlocal-to-local convergence belongs to a broader line of research on nonlocal-to-local limits, which has been studied, among others, for the Cahn--Hilliard equation \cite{MR4574535, MR4927649, MR4749398, MR4198717}, the porous medium and fast diffusion equations \cite{carrillo2024nonlocal, MR4732300, MR3913840, MR4858611}, hyperbolic conservation laws \cite{MR3961295, MR4553943}, and PDEs related to the machine learning algorithms \cite{carrillo2026stein,MR4646041}.\\

\noindent The only works addressing the convergence of SPH with a reasonably mild (in particular, algebraic) relationship between the nonlocality parameter $\eps$ and the number of particles $N$ are those of Oelschl{\"{a}}ger \cite{MR1120850} and Franz and Wendland \cite{MR3851044}. They established convergence of SPH approximations for quadratic pressure (i.e. $\gamma=2$), assuming that the limiting Euler system admits a smooth solution. The latter assumption is satisfied locally in time by classical Kato's theory \cite{MR390516}. The restriction to $\gamma=2$ is particularly important, as it makes the pressure force linear in the density in the velocity equation, allowing the use of a quadratic relative energy functional. More precisely, the convergence analysis relies on a relative entropy argument, based on controlling a quantity of the form
\[
Q(t)=\frac12\sum_{i=1}^N m_i|V_i(t)-u(t,X_i(t))|^2
+\frac12\|L_\eps\ast\rho_t^N-\rho(t)\|_{L^2}^2,
\]
where $\rho_t^N=\sum_{i=1}^N m_i\delta_{X_i(t)}$ is the empirical particle density and $(\rho, u)$ solves \eqref{sys:euler_general}. The key step is to derive a differential inequality for $Q(t)$, which yields quantitative convergence estimates through Grönwall's lemma and Fourier analysis. Moreover, Franz and Wendland \cite{MR4222186} showed that suitable choices of smoothing kernels, satisfying higher-order moment and approximation conditions, lead to improved convergence rates. Although we do not pursue this direction in the present work, it would be interesting to investigate whether analogous kernel constructions could also improve the convergence rates obtained here.  \\

\noindent In this paper, we bring together the analysis of nonlocal Euler equations and their particle approximations to provide a unified convergence framework. Our approach is based on the relative entropy method and applies to nonlinear pressure laws $\rho^\gamma$ with $\gamma\geq 2$. The relative entropy method is an important tool first introduced by Dafermos in the context of scalar conservation laws \cite{dafermos1979thesecond, dafermos2016hyperbolic} and helps describe how two physical systems differ in time when they start from close initial states. At this point it has found use in a various different fields, ranging from weak-strong uniqueness problems \cite{woznicki2022weak, wiedemann2017fluid, brenier2011weak, gwiazda2015weakstrong}, to stability studies, asymptotic limits, dimension reduction \cite{bella2014dimension, christoforou2018relative, feireisl2012relative, giesselmann2014singular, giesselemann2017stability, lattanzio2017fromgas, gallenmuller2024cahn}. An important advantage of the proposed nonlocal pressure operator is that it preserves the energy structure of the original Euler system, a property that is also essential to our analysis. For quadratic pressure ($\gamma=2$), the underlying argument is essentially that of Oelschl{\"{a}}ger \cite{MR1120850} and Franz and Wendland \cite{MR3851044,MR4222186}, which we formulate in a more general setting allowing for measure-valued densities in the nonlocal system. This encompasses both continuum solutions and particle approximations within the same framework. Our main contribution is to extend this analysis to nonlinear pressure laws with $\gamma>2$, establishing quantitative convergence rates towards smooth solutions of the local compressible Euler equations as $\eps\to0$, with algebraic, rather than exponential, dependence on the nonlocality parameter $\eps$ and the number of particles $N$ in the SPH setting. In addition, we obtain algebraic convergence rates for the density in the Wasserstein distance and negative Sobolev norms, as well as for the momentum in the flat metric. This extension presents significant analytical challenges, since, unlike in the quadratic case, the relative entropy cannot be estimated directly using Plancherel's theorem and Fourier analysis. Instead, our proof relies on careful decompositions of the nonlinear terms and estimates involving Riesz transforms. To the best of our knowledge, this provides the first rigorous convergence result for SPH in this nonlinear-pressure regime with an algebraic relationship between $\eps$ and $N$.\\

\noindent The remainder of the introduction is organized as follows. Section \ref{intro:sect_quadratic_pressure} presents the nonlocal-to-local convergence result for continuous solutions on the whole space $\R^d$ in the case $\gamma=2$ (Theorem \ref{thm:main}). Section \ref{sect_intro_general_pressure} discusses an extension of this result to $\gamma\geq 2$ on the torus $\T^d$ (Theorem \ref{thm:main_per}). Section \ref{sect_intro_sph_particles} establishes convergence rates for the SPH method (Theorem~\ref{thm:main_particles}) as a direct consequence of the quantitative estimates obtained in Theorems \ref{thm:main} and \ref{thm:main_per}. Finally, Section~\ref{sect:intro_notation} introduces the notation used throughout the paper. The remainder of the paper is organized as follows. Section \ref{sect:def_weak_slns_SPH_approx} introduces the notion of solutions and explains how the SPH approximation fits into this framework. Sections~\ref{sect:proof_theorem_quadratic_pressure}, \ref{sect:proof_theorem_general_pressure}, and \ref{sect:proof_main_result_for_particles} are devoted to the proofs of Theorems \ref{thm:main}, \ref{thm:main_per}, and \ref{thm:main_particles}, respectively. Finally, Appendix contains auxiliary results.

\subsection{The case of the quadratic pressure on $\R^d$}\label{intro:sect_quadratic_pressure}
We consider a system 
\begin{equation}\label{sys:non_loc_euler}
\begin{aligned}
    \partial_t\rho + \DIV_x(\rho\, u) &= 0 &&\quad\mbox{in }(0,T)\times \R^d,\\
    \partial_t (\rho\,u)  + \DIV_x (\rho\, u \otimes u ) + \rho\nabla_x (K_\varepsilon\star\rho) + \kappa\,\rho u  &= 0 &&\quad\mbox{in }(0,T)\times \R^d, \\
    \end{aligned}
\tag{NLE}
\end{equation}
where $\kappa \in \R$ and we study its convergence as $\varepsilon\to 0^+$ to
\begin{equation}\label{sys:loc_euler}
\begin{aligned}
    \partial_t\rho + \DIV_x(\rho\, u) &= 0 &&\quad\mbox{in }(0,T)\times \R^d,\\
    \partial_t (\rho\,u)  + \DIV_x (\rho\, u \otimes u ) + \frac{1}{2}\nabla_x(\rho^2) + \kappa\,\rho u  &= 0 &&\quad\mbox{in }(0,T)\times \R^d. \\
    \end{aligned}
\tag{E}
\end{equation}
The concept of solutions is deferred to Definition~\ref{def:weak_solution} below. Here, the kernel $K_{\eps}(x) = \frac{1}{\varepsilon^d}\, K(x/\varepsilon)$, where we assume that the kernel $K$ is given by
\begin{equation}\label{eq:particular_form_of_the_kernel_quadratic_press_part1}
    K(x) = L_m\star L_m(x),\quad m > \max\{d, 2\},
\end{equation}
where the kernel $L_m$ is the so-called Bessel's potential of order $m$, that is $L_m\in L^1(\R^d)$ is such that its Fourier transform satisfies
\begin{equation}\label{eq:particular_form_of_the_kernel_quadratic_press_part2}
    \hat{L}_m(\xi) := \frac{1}{(1 + 4\pi^2|\xi|^2)^{m/2}}.
\end{equation}
It is clear that $L_m(-x) = L_m(x)$. Denoting $L_{m, \varepsilon}(x) = \frac{1}{\varepsilon^d}L_m(x/\varepsilon)$, we have
    $$
    K_\varepsilon := L_{m, \varepsilon}\star L_{m, \varepsilon}, \qquad
    \hat{K}(\xi) = (\hat{L}_m(\xi))^2,\qquad \hat{K}_\varepsilon(\xi) = \hat{K}(\varepsilon\, \xi) = \frac{1}{(1 + 4\pi^2|\varepsilon\,\xi|^2)^m}, \qquad K_{\eps}(-x) = K_{\eps}(x).
    $$
Our main result in this case reads:
\begin{thm}\label{thm:main}
    Let $\kappa \in \R$. Suppose that $(\rho^\varepsilon, u^\varepsilon, \mu^\varepsilon)$ is a measure solution to the system \eqref{sys:non_loc_euler} as in Definition~\ref{def:weak_solution} with initial datum $(\rho^\varepsilon_0, u^\varepsilon_0)$, and that $(r, U)$ is a smooth solution to the system \eqref{sys:loc_euler} with smooth initial conditions $r_0\in L^1(\R^d)$, $U_0\in L^2(\R^d, \diff r_0)$ satisfying
    $$
    \int_{\R^d}|x|^2\, r_0\diff x < +\infty, \qquad \int_{\R^d}r_0\diff x = 1,
    $$
    such that $r > 0$, $\nabla_x^2 r\in L^\infty((0, T)\times \R^d)$, $U\in L^\infty(0, T; W^{2,\infty}(\R^d))$. Define the relative entropy
    \begin{align*}
        \mathcal{E}_{rel}^\varepsilon(t) := \frac{1}{2}\int_{\R^d}|u^\varepsilon - U|^2\diff\rho^\varepsilon_t(x) + \frac{1}{2}\int_{\R^d}|L_{m, \varepsilon}\star(\diff\rho^\varepsilon - \diff r)|^2\diff x.
    \end{align*}
Then, there exist constants $C_1, C_2 > 0$ depending on the norms of the strong solution $(r, U)$, dimension $d$, exponent $m$, time interval $T$ (the constant $C_2$ additionally depends on the $\|u_0^\eps\|_{L^2(\R^d;\diff\rho_0^\eps)}$) such that for a.e. $t\in (0, T)$
    \begin{equation}\label{eq:energy_relative_main_bound_Rd}
        \mathcal{E}_{rel}^\varepsilon(t) + \frac{1}{2}\mathrm{tr}(\mu^\varepsilon_t(\R^d)) + \kappa \int_{0}^t\int_{\R^d}|u^\varepsilon - U|^2\diff \rho^\varepsilon_\tau(x)\diff\tau \leq \mathcal{E}_{rel}^\varepsilon(0) + C_1\int_0^t\mathcal{E}_{rel}^\varepsilon(\tau) + + \frac{1}{2}\mathrm{tr}(\mu^\varepsilon_\tau(\R^d))\diff\tau + C_2\, \varepsilon,
    \end{equation}
\noindent \underline{Case $\kappa \geq 0$.}  There is a constant $C_3=C\left(
T,C_1,C_2,
\|U\|_{L^\infty(0,T;W^{1,\infty}(\R^d))}
\right)$ and a numerical constant $C_4$, such that
    \begin{equation}\label{eq:energy_relative_main_bound_Rd_gronwall_estimate}
\mathcal{E}_{\mathrm{rel}}^\varepsilon(t) +\frac{1}{2}\mathrm{tr}(\mu^\varepsilon_t(\R^d))
+\kappa \int_0^t\int_{\mathbb{R}^d}
|u^\varepsilon-U|^2\,\diff\rho^\varepsilon_\tau(x)\diff\tau
\leq
e^{C_1 t}\left(
\mathcal{E}_{\mathrm{rel}}^\varepsilon(0)+C_2\, \varepsilon
\right),
\end{equation}
    \begin{equation}\label{eq:final_estimate_W_2_on_R_d}
    W_2^2(\rho_t^\varepsilon,r_t)
\leq
C_3\left(
W_2^2(\rho_0^\varepsilon,r_0)
+\mathcal E_{\mathrm{rel}}^\varepsilon(0)
+\varepsilon
\right),
    \end{equation} 
    \begin{equation}\label{eq:final_estimate_flat_norm_on_R_d}
    d_f^2(\rho^\varepsilon_t\, u^\varepsilon_t, r_t\,U_t) \leq 
C_3\left(
W_2^2(\rho_0^\varepsilon,r_0)
+\mathcal E_{\mathrm{rel}}^\varepsilon(0)
+\varepsilon
\right),
    \end{equation}
    \begin{equation}\label{eq:final_estimate_neg_SS_on_R_d}
    \|\rho^\varepsilon_t - r_t\|^2_{H^{-m}(\R^d)} + \|\mu^\varepsilon_t\|_{TV} \leq C_4\, e^{C_1 t}\left(
\mathcal{E}_{\mathrm{rel}}^\varepsilon(0)+C_2\, \varepsilon
\right),
    \end{equation} 
    where $W_2$ is the 2-Wasserstein metric and $d_f$ is the flat metric, see Appendix~\ref{sect:wasserstein_flat}. Assume further, that the initial conditions $(\rho_0^\varepsilon, u_0^\varepsilon)$ are well-prepared in the sense that $\mathcal{E}_{rel}^\varepsilon(0)\to 0$ as $\varepsilon \to 0$ and $\|u_0^\eps\|_{L^2(\R^d;\diff\rho_0^\eps)}$ is bounded independently of $\varepsilon$. Then, $\mathcal{E}_{rel}^\varepsilon(t)\to 0$ and $\rho^{\varepsilon}_t \to r_t$ in $H^{-m}(\R^d)$ for a.e. $t \in (0,T)$. Additionally, if $W_2(\rho_0^{\varepsilon},r_0) \to 0$, then we have $W_2^2(\rho_t^\varepsilon,r_t) \to 0$ and $d_f^2(\rho^\varepsilon_t\, u^\varepsilon_t, r_t\,U_t) \to 0$.\\
\underline{Case $\kappa < 0$.} The same conclusion is true but the inequality \eqref{eq:energy_relative_main_bound_Rd_gronwall_estimate} read
   \begin{equation}\label{eq:energy_relative_main_bound_Rd_gronwall_estimate_negative_k}
\mathcal{E}_{\mathrm{rel}}^\varepsilon(t) + \frac{1}{2}\mathrm{tr}(\mu^\varepsilon_t(\R^d)) \leq
e^{(C_1 - 2\kappa) t}\left(
\mathcal{E}_{\mathrm{rel}}^\varepsilon(0)+C_2\, \varepsilon
\right),
\end{equation}
Consequently, \eqref{eq:final_estimate_W_2_on_R_d}--\eqref{eq:final_estimate_neg_SS_on_R_d} hold true with constant $C_1 := C_1 -2\kappa$ and $C_3$ depends additionally on $\kappa$. 
\end{thm}

\noindent For $\gamma=2$, the nonlocal-to-local convergence for sufficiently regular initial densities was already established by Di~Lisio, Grenier, and Pulvirenti \cite{MR1487955} and Danchin and Mucha \cite{MR4764519}, by deriving uniform-in-$\eps$ estimates in Sobolev and Besov spaces, respectively. For measure-valued densities, convergence was established by Oelschl\"{a}ger \cite{MR1120850} and subsequently studied by Franz and Wendland \cite{MR3851044,MR4222186}. Our proof essentially follows the latter approach. We express the relative entropy in Fourier variables and estimate the resulting terms using Taylor expansions, which rely on the specific form of the smoothing kernel \eqref{eq:particular_form_of_the_kernel_quadratic_press_part1}--\eqref{eq:particular_form_of_the_kernel_quadratic_press_part2}. Although the argument could potentially be extended to a broader class of kernels, we regard the present result primarily as a proof of concept. We also note that a similar approach based on Taylor expansions has been employed in the study of the nonlocal-to-local limit for Stein variational gradient descent (SVGD) \cite{carrillo2026stein}.

\subsection{The case of the pressure $\rho^{\gamma}$, $\gamma \geq 2$}\label{sect_intro_general_pressure} We define $\T^d := (0, 1)^d|_{\{0, 1\}}$ and we consider the periodization of the kernel $L_m$ given by
\begin{equation}\label{eq:nonlocal_kernel_periodic}
    L^{per}_m := \sum_{l\in\I^d}L_m(x +l),\qquad L^{per}_{m,\eps}(x) := \sum_{l\in\I^d}L_{m,\eps}(x +l) \qquad x\in\T^d.
\end{equation}
The reason is that the function $L_m$ is globally supported, hence even for small $\eps>0$, $L_{m,\eps}$ will not be interpretable as a periodic function. We consider a system 
\begin{equation}\label{sys:non_loc_euler_per}
\begin{aligned}
    \partial_t\rho + \DIV_x(\rho\, u) &= 0 &&\quad\mbox{in }(0,T)\times \T^d,\\
    \partial_t (\rho\,u)  + \DIV_x (\rho\, u \otimes u ) + \rho\nabla_x (L^{per}_{m,\eps}\star(L_{m,\eps}^{per}\star\rho)^{\gamma-1}) + \kappa \,\rho\, u  &= 0 &&\quad\mbox{in }(0,T)\times \T^d, \\
    \end{aligned}
\tag{NLE-P}
\end{equation}
and its convergence as $\varepsilon\to 0^+$ to
\begin{equation}\label{sys:loc_euler_per}
\begin{aligned}
    \partial_t\rho + \DIV_x(\rho\, u) &= 0 &&\quad\mbox{in }(0,T)\times \T^d,\\
    \partial_t (\rho\,u)  + \DIV_x (\rho\, u \otimes u ) + \frac{1}{\gamma}\nabla_x\rho^\gamma + \kappa \,\rho\, u  &= 0 &&\quad\mbox{in }(0,T)\times \T^d. \\
    \end{aligned}
\tag{E-P}
\end{equation}
\noindent Our main result in this case reads:
\begin{thm}\label{thm:main_per}
    Let $\gamma \geq 2$ and $\kappa \in \R$. Suppose that $(\rho^\varepsilon, u^\varepsilon, \mu^\varepsilon)$ is a measure solution to the system \eqref{sys:non_loc_euler_per} as in Definition \ref{def:weak_solution_per} with initial datum $(\rho^\varepsilon_0, u^\varepsilon_0)$, and that $(r, U)$ is a smooth solution to the system \eqref{sys:loc_euler_per} with smooth initial conditions $r_0\in L^1(\T^d)$, $U_0 \in L^2(\T^d, \diff r_0)$ satisfying
    $$
    r_0\geq r_* > 0, \qquad \int_{\T^d}r_0\diff x = 1,
    $$
    such that $r\in L^\infty(0, T; W^{2,\infty}(\T^d))$, $U\in L^\infty(0, T; W^{2,\infty}(\T^d))$. Define for $f, g \geq 0$, $\beta > 1$
    \begin{align}\label{eq:gamma_distance}
d_\beta(f, g) = f^\beta - \beta\, g^{\beta - 1}f + (\beta - 1)\,g^\beta,
\end{align}
and the relative entropy
    \begin{align*}
        \mathcal{E}_{rel}^\varepsilon(t) := \frac{1}{2}\int_{\T^d}|u^\varepsilon - U|^2\diff\rho^\varepsilon_t(x) + \frac{1}{\gamma}\int_{\T^d}d_\gamma(L^{per}_{m, \varepsilon}\star\diff\rho^\varepsilon, L_{m,\eps}^{per} \star\diff r)\diff x.
    \end{align*}
    Then, there are constants $C_1, C_2 > 0$ depending on the norms of the strong solution $(r, U)$, dimension $d$, exponents $\gamma$ and $m$, time interval $T$ (the constant $C_2$ additionally depends on the $\|u_0^\eps\|_{L^2(\T^d;\diff\rho_0^\eps)}$) such that for a.e. $t\in(0,T)$
    \begin{equation}\label{eq:rel_energy_before_Gronwall_general_pressure}
        \mathcal{E}_{rel}^\varepsilon(t) + \frac{1}{2}\mathrm{tr}(\mu^\varepsilon_t(\T^d)) + \kappa \int_{0}^t\int_{\T^d}|u^\varepsilon - U|^2\diff \rho^\varepsilon_\tau(x)\diff\tau \leq \mathcal{E}_{rel}^\varepsilon(0) + C_1\int_0^t\mathcal{E}_{rel}^\varepsilon(\tau) + \frac{1}{2}\mathrm{tr}(\mu^\varepsilon_\tau(\T^d))\diff\tau + C_2\, \eps.
    \end{equation}
\underline{Case $\kappa \geq 0$.} There are constants $C_3=C\left(
T,C_1,C_2,
\|U\|_{L^\infty(0,T;W^{1,\infty}(\T^d))}
\right)$, $C_4 = C_4(r_*, \|r\|_{\infty})$, and a numerical constant $C_5$, such that
\begin{equation}\label{eq:relative_entropy_after_gronwall_nonlinear_press_k_pos}
\mathcal{E}_{\mathrm{rel}}^\varepsilon(t) + \frac{1}{2}\mathrm{tr}(\mu^\varepsilon_t(\T^d))
+\kappa \int_0^t\int_{\mathbb{T}^d}
|u^\varepsilon-U|^2\,\diff\rho^\varepsilon_\tau(x)\diff\tau
\leq
e^{C_1 t}\left(
\mathcal{E}_{\mathrm{rel}}^\varepsilon(0)+C_2\, \eps
\right).
\end{equation} 
\begin{equation}\label{eq:wasserstein_rate_nonlinear_press}
    W_2^2(\rho_t^\varepsilon,r_t)
\leq
C_3\left(
W_2^2(\rho_0^\varepsilon,r_0)
+\mathcal E_{\mathrm{rel}}^\varepsilon(0)
+ \eps
\right),
\end{equation} 
\begin{equation}\label{eq:flat_rate_nonlinear_press}
    d_f^2(\rho^\varepsilon_t\, u^\varepsilon_t, r_t\,U_t) \leq 
C_3\left(
W_2^2(\rho_0^\varepsilon,r_0)
+\mathcal E_{\mathrm{rel}}^\varepsilon(0)
+  \eps
\right),
\end{equation}
\begin{equation}\label{eq:negSS_rate_nonlinear_press}
    \|\rho^\varepsilon_t - r_t\|^2_{H^{-m}(\R^d)} \leq C_4\, e^{C_1 t}\left(
\mathcal{E}_{\mathrm{rel}}^\varepsilon(0)+C_2\, \eps \right), \qquad   \|\mu^\varepsilon_t\|_{TV} \leq C_5\, e^{C_1 t}\left(
\mathcal{E}_{\mathrm{rel}}^\varepsilon(0)+C_2\, \eps\right),
\end{equation}
where $W_2$ is the 2-Wasserstein metric and $d_f$ is the flat metric, see Appendix~\ref{sect:wasserstein_flat}. Assume further, that the initial conditions $(\rho_0^\varepsilon, u_0^\varepsilon)$ are well-prepared in the sense that $\mathcal{E}_{rel}^\varepsilon(0)\to 0$ as $\varepsilon \to 0$ and $\|u_0^\eps\|_{L^2(\T^d;\diff\rho_0^\eps)}$ is bounded independently of $\varepsilon$. Then, $\mathcal{E}_{rel}^\varepsilon(t)\to 0$ and $\rho^{\varepsilon}_t \to r_t$ in $H^{-m}(\R^d)$ for a.e. $t \in (0,T)$. Additionally, if $W_2(\rho_0^{\varepsilon},r_0) \to 0$, then we have $W_2^2(\rho_t^\varepsilon,r_t) \to 0$ and $d_f^2(\rho^\varepsilon_t\, u^\varepsilon_t, r_t\,U_t) \to 0$.\\
\underline{Case $\kappa < 0$.} The same conclusion is true but the inequality \eqref{eq:relative_entropy_after_gronwall_nonlinear_press_k_pos} read
   \begin{equation}\label{eq:relative_entropy_after_gronwall_nonlinear_press_k_neg}
\mathcal{E}_{\mathrm{rel}}^\varepsilon(t) + \frac{1}{2}\mathrm{tr}(\mu^\varepsilon_t(\T^d)) \leq
e^{(C_1-2\kappa) t}\left(
\mathcal{E}_{\mathrm{rel}}^\varepsilon(0)+C_2\,\eps
\right).
\end{equation}
Consequently, \eqref{eq:wasserstein_rate_nonlinear_press}--\eqref{eq:negSS_rate_nonlinear_press} hold true with constant $C_1 := C_1 -2\kappa$ and $C_3$ depends additionally on $\kappa$. 
\end{thm}
\noindent For nonlinear pressure laws with $\gamma>2$, the only previously known nonlocal-to-local convergence result, to the best of our knowledge, is due to Di Lisio, Grenier, and Pulvirenti \cite{MR1487955}, who obtained uniform-in-$\eps$ Sobolev estimates for sufficiently regular solutions. Our proof instead relies on the relative entropy method, combined with a careful decomposition of the nonlinear terms and estimates based on Fourier analysis and Riesz transforms. A crucial restriction of our approach is that it requires $\gamma\geq 2$ and a bounded spatial domain. Indeed, a key ingredient is the inequality
$
|\rho-r|^2\leq C(r_*, \gamma)\,d_\gamma(\rho,r),
$
where $r_*>0$ is a lower bound on the density $r$, which is true for all $\gamma \geq 2$, $\rho\geq0$ and $r \geq r_*$. The restriction of a bounded domain allows us to assume a strictly positive lower bound on the limiting density without imposing infinite total mass.

\subsection{Smoothed particle hydrodynamics}\label{sect_intro_sph_particles} We now consider a measure $\rho = \sum_{i=1}^N m_i\, \delta_{X_i(t)}$ which will be a~solution to the nonlocal Euler equation. On $\R^d$, in order to solve \eqref{sys:non_loc_euler}, we define
\begin{equation}\label{eq:SPH_Rd}
    \begin{split}
        \dot{X_i}(t) &= V_i(t),\\
        \dot{V_i}(t) &=  - \sum_{j=1}^N m_j\, \nabla K_{\eps}(X_i(t)-X_j(t)) - \kappa \,V_i(t).
    \end{split}
\end{equation}
On the other hand, on $\T^d$, in order to solve \eqref{sys:non_loc_euler_per}, we define 
\begin{equation}\label{eq:SPH_Td}
    \begin{split}
        \dot{X_i}(t) &= V_i(t),\\
        \dot{V_i}(t) &= -\int_{\T^d} \nabla_x L^{per}_{m,\eps}(X_i(t) - y) \, \left(\sum_{k=1}^N m_k\, L_{m,\eps}^{per}(y - X_k(t)) \right)^{\gamma-1}\diff y - \kappa \,V_i(t). 
    \end{split}
\end{equation}
It is proved in Proposition~\ref{prop:particle_sln_is_a_sol} that $(\rho, u)$ with $u(t,X_i(t)) := V_i(t)$ and $(X_i(t), V_i(t))$ given by \eqref{eq:SPH_Rd} is a measure solution to \eqref{sys:non_loc_euler} (and similarly, if $(X_i(t), V_i(t))$ given by \eqref{eq:SPH_Td}, then $(\rho, u)$ solves \eqref{sys:non_loc_euler_per}). Moreover, there is no concentration component i.e. $\mu = 0$. Therefore, one immediately obtains the following result.

\begin{thm}\label{thm:main_particles} 
{\bf (Case $\R^d$, $\gamma = 2$)} Let $r^0$, $U^0$ be smooth initial conditions with $r^0 \in \mathcal{P}(\R^d)$, $U^0 \in L^2(\R^d, \diff r^0) \cap L^{\infty}(\R^d)$, and $\int_{\R^d} |x|^{2+\delta} \diff r(x) < \infty$ for some $\delta>0$. Let $m>d$. Let $(r,U)$ be the smooth solution of \eqref{sys:loc_euler} with initial conditions $r^0, U^0$ such that $r > 0$, $\nabla_x^2 r\in L^\infty((0, T)\times \R^d)$, $U\in L^\infty(0, T; W^{2,\infty}(\R^d))$. Let $r^{N}_0=\sum_{i=1}^{k^d+1} r(A_i)\, \delta_{X_i(0)}$ be the approximation of $r^0$ given by Lemma \ref{lem:particle_approx_of_a_gen_measure} and let $N := k^d+1$. Consider $(X_i, V_i)_{i=1}^N \subset \R^{2N}$ be the solution to \eqref{eq:SPH_Rd} with initial conditions $(X_i(0), V_i(0))_{i=1}^N$ satisfying $V_i(0)= U^0(X_i(0))$. Let $\rho^{N,\varepsilon}_t = \sum_{i=1}^{k^d+1} r(A_i)\, \delta_{X_i(t)}$ and let $u^N_t$ be the function defined on the support of $\rho^{N,\varepsilon}_t$ via $u^{N,\varepsilon}_t(X_i(t)) := V_i(t)$. Then, there is a constant $C_{tot}$ depending on the norms of the strong solution $(r, U)$, dimension $d$, exponent $m$, coefficient $\kappa$, time interval $T$ such that 
    \begin{equation}\label{eq:rate_SPH_Rd_W2}
    W_2^2(\rho_t^{N,\varepsilon},r_t)
\leq
C_{tot} \left(
W_2^2(\rho_0^{N, \varepsilon},r_0)
+ {\varepsilon^{-(2+d)}}\, k^{-\frac{2\,\delta}{(2+\delta)}}
+\varepsilon
\right),\qquad N = k^d + 1
    \end{equation}
    \begin{equation}\label{eq:rate_SPH_Rd_fn}
    d_f^2(\rho^{N, \varepsilon}_t\, u^{N, \varepsilon}_t, r_t\,U_t) \leq 
C_{tot} \left(
W_2^2(\rho_0^{N, \varepsilon},r_0)
+{\varepsilon^{-(2+d)}}\, k^{-\frac{2\,\delta}{(2+\delta)}}
+\varepsilon
\right),\qquad N = k^d + 1
    \end{equation} 
    \begin{equation}\label{eq:rate_SPH_Rd_negSS}
    \|\rho^{N,\varepsilon}_t - r_t\|^2_{H^{-m}(\R^d)} \leq C_{tot}\left(
{\varepsilon^{-(2+d)}}\, k^{-\frac{2\,\delta}{(2+\delta)}}+\varepsilon
\right),\qquad N = k^d + 1.
    \end{equation}
{\bf (Case $\T^d$, $\gamma \geq 2$)} Let $\gamma \geq 2$, $m\geq d+1$ and let $r^0$, $U^0$ be smooth initial conditions with $r^0 \in \mathcal{P}(\T^d)$, $U^0 \in L^2(\T^d, \diff r^0) \cap L^{\infty}(\T^d)$, $r_0 \geq r_{*}$. Let $r\in L^\infty(0, T; W^{2,\infty}(\T^d))$, $U\in L^\infty(0, T; W^{2,\infty}(\T^d))$ be the corresponding solution to \eqref{sys:loc_euler_per}. Let $r^{N}_0=\sum_{i=1}^{k^d} r(A_i)\, \delta_{X_i(0)}$ be the approximation of $r^0$ given by Lemma \ref{lem:particle_approx_of_a_gen_measure} and let $N := k^d$. Consider $(X_i, V_i)_{i=1}^N \subset \R^{2N}$ be the solution to \eqref{eq:SPH_Td} with initial conditions $(X_i(0), V_i(0))_{i=1}^N$ satisfying $V_i(0)= U^0(X_i(0))$. Let $\rho^{N,\varepsilon}_t = \sum_{i=1}^{k^d} r(A_i)\, \delta_{X_i(t)}$ and let $u^{N,\varepsilon}_t$ be the function defined on the support of $\rho^{N,\varepsilon}_t$ via $u^{N,\varepsilon}_t(X_i(t)) := V_i(t)$.  Then, there is a constant $C_{tot}$ depending on the norms of the strong solution $(r, U)$, dimension $d$, exponents $\gamma$ and $m$, coefficient $\kappa$, lower bound $r_*$, and time interval $T$ such that 
\begin{equation}\label{eq:rate_SPH_Td_W2}
    W_2^2(\rho_t^{N, \varepsilon},r_t)
\leq
C_{tot} \left(
k^{-2}
+k^{-1}\, \eps^{-(2d+1)}
+ \eps
\right),\qquad N = k^d.
\end{equation}
\begin{equation}\label{eq:rate_SPH_Td_fn}
    d_f^2(\rho^{N,\varepsilon}_t\, u^{N,\varepsilon}_t, r_t\,U_t) \leq 
C_{tot} \left(
k^{-2}
+k^{-1}\,  \eps^{-(2d+1)}
+  \eps
\right),\qquad N= k^d.
\end{equation}
\begin{equation}\label{eq:rate_SPH_Td_negSS}
    \|\rho^{N,\varepsilon}_t - r_t\|^2_{H^{-m}(\R^d)}  \leq C_{tot} \left(
k^{-1}\,  \eps^{-(2d+1)}+ \eps
\right),\qquad N = k^d.
\end{equation}

\end{thm}

\noindent To the best of our knowledge, the only previous rigorous convergence results for SPH with an algebraic relationship between the number of particles $N$ and the smoothing parameter $\eps$ are those of Oelschl{\"{a}}ger \cite{MR1120850} and Franz and Wendland \cite{MR3851044,MR4222186}, which are restricted to quadratic pressure ($\gamma=2$). A notable feature of our approach is that the convergence results for the particle system follow directly from our general estimates for measure-valued solutions of the nonlocal Euler equations. Thus, no separate analysis of the particle dynamics is required, and the same framework yields quantitative convergence rates for both continuum and particle approximations, including rates for the density in the Wasserstein distance and negative Sobolev norms, and for the momentum in the flat norm. In particular, our results for nonlinear pressure laws with $\gamma>2$ are entirely new and, to the best of our knowledge, provide the first rigorous convergence results for SPH in this regime with an algebraic relationship between $N$ (via $N=k^d+1$) and $\eps$.

\subsection{Notation}\label{sect:intro_notation} Throughout the article we will use $x\in \R^d$ or $x\in\T^d$ to denote a spatial variable for a dimension $d\geq 1$, and $t\in (0, T)$ for a time variable with time horizon $T > 0$. For any two vectors $a,b\in \R^d$, by $a\cdot b$ we will mean their standard scalar product i.e. $a\cdot b := \sum_{i = 1}^d a_i\, b_i$, and for any two matrices $A, B\in \R^{d\times d}$, by $A: B$ we will denote the scalar product emanating from the Frobenius norm, that is $A: B := \mathrm{tr}(A^T B) = \sum_{i, j = 1}^d A_{ij}\, B_{ij}$. When it comes to functional spaces we employ the standard notation for Lebesgue's ($L^p$) and Sobolev spaces ($W^{k, p}$). Furthermore, we use $\mathcal{M}$ to represent the space of Radon measures and $\mathcal{P}$ for its subset of probability measures. We will frequently utilize $\hat{f}$ to signify the Fourier transform of a function or a measure, i.e.
$$
\hat{f}(\xi) := \int_{\R^d}f(x)\exp(-2\pi x\cdot \xi)\diff x.
$$
We also need a definition of symmetric, matrix-valued measures. The main property of such measures is that their total variation is controlled by their trace.
\begin{Def}
    Suppose $\Omega = \R^d$ or $\T^d$. We say that $\mu\in \mathcal{M}^+(\Omega; \R^{d\times d}_{\mathrm{sym}})$, iff
    $$
    \int_{\Omega}\phi(x)\xi\otimes\xi:\diff\mu \geq 0,
    $$
    for any $\xi\in \R^d$, $\phi\in C_c(\Omega)$, $\phi\geq 0$.
\end{Def}
\begin{prop}[Lemma 3.2, \cite{woznicki2022weak}]\label{prop:bound_symmetric_measure_trace}
    Suppose that $\mu\in \mathcal{M}^+(\Omega; \R^{d\times d}_{\mathrm{sym}})$ for $\Omega = \R^d$ or $\T^d$, then there exists a constant $C > 0$ (independent of $\mu$) such that
    $$
    |\mu|(\Omega) \leq C\,\mathrm{tr}\left(\mu(\Omega)\right).
    $$
\end{prop}

\section{Notions of solutions and the SPH approximation}\label{sect:def_weak_slns_SPH_approx}
\noindent We define here measure solution to the systems \eqref{sys:non_loc_euler} and \eqref{sys:non_loc_euler_per}. This are motivated by two facts. First, the measure $\mu$ is included in the momentum equation to allow for possible oscillations when one constructs solutions via an approximating sequence. Indeed, the only global solutions to \eqref{sys:non_loc_euler} and \eqref{sys:non_loc_euler_per} are known to exist only at the level of measure-valued solutions where additional terms present in the formulations account for possible oscillations of the approximating sequences. Second, we would like to include the particles solutions to these equations, which is reflected by the fact that density $\rho_t$ is required to be only a measure for a.e. $t \in (0,T)$.  
\begin{Def}\label{def:weak_solution}
We say that $(\rho, u, \mu)$ is a measure solution to the system \eqref{sys:non_loc_euler} with initial datum $(\rho_0, u_0)$ if
    $$
    \rho\in L^\infty(0, T; \mathcal{P}(\R^d)),\quad u\in L^\infty(0, T; L^2(\R^d;\diff\rho_t)),\quad \mu\in L^\infty(0, T; \mathcal{M}^+(\R^d; \R^{d\times d}_{\mathrm{sym}}))
    $$
    and the following are satisfied:
    \begin{itemize}
        \item continuity equation: for any $\varphi\in C^\infty_c([0, T)\times \R^d)$ and a.e. $t\in (0, T)$
        \begin{align}\label{eq:def_weak_cont}
            \int_{\R^d}\varphi(t, x)\diff\rho_t(x)  - \int_{\R^d}\varphi(0, x)\diff\rho_0(x) = \int_0^t\int_{\R^d} \big(\partial_t\varphi + u\cdot \nabla_x\varphi \big) \diff \rho_\tau(x)\diff\tau,
        \end{align}
        \item momentum equation: for any $\phi\in C^\infty_c([0, T)\times \R^d)$ and a.e. $t\in (0, T)$
        \begin{equation}\label{eq:def_weak_moment}
            \begin{split}
                &\int_{\R^d}\phi(t, x)\cdot u(t,x)\diff\rho_t(x) - \int_{\R^d}\phi(0, x)\cdot u_0(x)\diff\rho_0(x)\\
                &\qquad= \int_0^t\int_{\R^d} \big( u\cdot\partial_t\phi + u\otimes u:\nabla_x\phi - \kappa\,u\cdot\phi - \nabla_x(K_\varepsilon\star\diff\rho)\cdot\phi \big) \diff\rho_\tau(x)\diff \tau + \int_0^t\int_{\R^d}\nabla_x\phi : \diff\mu_\tau(x)\diff \tau,
            \end{split}
        \end{equation}
        \item energy inequality: for a.e. $t\in (0, T)$
        \begin{equation}\label{ineq:def_weak_energy}
            \begin{split}
                \int_{\R^d}|u|^2(t, x)\diff\rho_t(x)& + \int_{\R^d}|L_{m, \varepsilon}\star\diff\rho|^2(t, x)\diff x + \mathrm{tr}(\mu_t(\R^d))\\
                & + 2\, \kappa \int_0^t\int_{\R^d}|u|^2\diff\rho_\tau(x)\diff\tau \leq \int_{\R^d}|u_0|^2\diff\rho_0(x) +  \int_{\R^d}|L_{m, \varepsilon}\star\diff\rho_0|^2(t, x)\diff x .
            \end{split}
        \end{equation}
    \end{itemize}
\end{Def}

\begin{Def}\label{def:weak_solution_per}
    We say that $(\rho, u, \mu)$ is a measure solution to the system \eqref{sys:non_loc_euler_per} with initial datum $(\rho_0, u_0)$ and kernel $L^{per}_m$ as in \eqref{eq:nonlocal_kernel_periodic} if
    $$
    \rho\in L^\infty(0, T; \mathcal{P}(\T^d)),\quad u\in L^\infty(0, T; L^2(\T^d;\diff\rho_t)),\quad \mu\in L^\infty(0, T; \mathcal{M}^+(\T^d; \R^{d\times d}_{\mathrm{sym}})),
    $$
    and the following are satisfied:
    \begin{itemize}
        \item continuity equation: for any $\varphi\in C^\infty_c([0, T)\times \T^d)$ and a.e. $t\in (0, T)$
        \begin{align}\label{eq:def_weak_cont_per}
            \int_{\T^d}\varphi(t, x)\diff\rho_t(x)  - \int_{\T^d}\varphi(0, x)\diff\rho_0(x)  = \int_0^t\int_{\T^d}\big(\partial_t\varphi + u\cdot \nabla_x\varphi \big) \diff \rho_\tau(x)\diff\tau,
        \end{align}
        \item momentum equation: for any $\phi\in C^\infty_c([0, T)\times \T^d)$ and a.e. $t\in (0, T)$
        \begin{equation}\label{eq:def_weak_moment_per}
            \begin{split}
                &\int_{\T^d}\phi(t, x)\cdot u(t,x)\diff\rho_t(x) - \int_{\T^d}\phi(0, x)\cdot u_0(x)\diff\rho_0(x) = \int_0^t\int_{\T^d} \big(u\cdot\partial_t\phi + u\otimes u:\nabla_x\phi\big) \diff\rho_\tau(x)\diff \tau \\
                & \qquad - \int_0^t\int_{\T^d}   \big(\kappa \,u\cdot\phi + \nabla_x(L^{per}_{m,\eps}\star(L^{per}_{m,\varepsilon}\star\diff\rho)^{\gamma - 1})\cdot\phi\big) \diff\rho_\tau(x)\diff \tau + \int_0^t\int_{\T^d}\nabla_x\phi : \diff\mu_\tau(x)\diff \tau,
            \end{split}
        \end{equation}
        \item energy inequality: for a.e. $t\in (0, T)$ \begin{equation}\label{ineq:def_weak_energy_per}
            \begin{split}
                \frac{1}{2}\int_{\T^d}|u|^2(t, x)\diff\rho_t(x)& + \frac{1}{\gamma}\int_{\T^d}(L^{per}_{m, \varepsilon}\star\diff\rho)^\gamma(t, x)\diff x +  \frac{1}{2}\mathrm{tr}(\mu_t(\T^d))\\
                & + \kappa\int_0^t\int_{\T^d}|u|^2\diff\rho_\tau(x)\diff\tau \leq \frac{1}{2}\int_{\T^d}|u_0|^2\diff\rho_0(x) +  \frac{1}{\gamma}\int_{\T^d}(L^{per}_{m, \varepsilon}\star\diff\rho_0)^\gamma(t, x)\diff x.
            \end{split}
        \end{equation}
    \end{itemize}
\end{Def}

\begin{prop}\label{prop:particle_sln_is_a_sol}
Let $N \geq 2$ and let $(X_i, V_i)_{i=1}^N$ be the solution to \eqref{eq:SPH_Rd} with some initial conditions $(X_i(0), V_i(0))_{i=1}^N$. Let $\rho = \sum_{i=1}^N m_i\, \delta_{X_i(t)}$ and let $u(t,X_i(t)) := V_i(t)$. Then, $(\rho, u, 0)$ is a measure solution to \eqref{sys:non_loc_euler} as in Definition~\ref{def:weak_solution}. Morover, it satisfies equality in \eqref{ineq:def_weak_energy} rather than only an inequality. Similarly, if \eqref{eq:SPH_Rd} is replaced by \eqref{eq:SPH_Td}, then $(\rho, u, 0)$ is a measure solution to \eqref{sys:non_loc_euler} as in Definition~\ref{def:weak_solution_per} and it satisfies an equality in \eqref{ineq:def_weak_energy_per}. 
\end{prop}
\begin{proof}
We compute
\begin{multline*}
\int_0^t\int_{\R^d}\big(\partial_t\varphi + u\cdot \nabla_x\varphi \big) \diff \rho_\tau(x)\diff\tau = \sum_{i=1}^N m_i \int_0^t \big( (\partial_t \varphi)(\tau, X_i(\tau)) + V_i(t) \cdot \nabla_x \varphi(\tau, X_i(\tau)) \big) \diff \tau = \\
 = \sum_{i=1}^N m_i \int_0^t \frac{\diff}{\diff \tau} \varphi(\tau, X_i(\tau)) \diff \tau =  \int_{\R^d}\varphi(t, x)\diff\rho_t(x)  - \int_{\R^d}\varphi(0, x)\diff\rho_0(x),
\end{multline*}
which yields \eqref{eq:def_weak_cont}. Similarly,
\begin{align*}
&\int_0^t\int_{\R^d} \big(u\cdot\partial_t\phi + u\otimes u:\nabla_x\phi\big) \diff\rho_\tau(x) \diff \tau  = \sum_{i=1}^N m_i \int_0^t \big(V_i(\tau) \cdot (\partial_t \phi)(\tau, X_i(\tau)) + V_i(\tau) \otimes V_i(\tau) : \nabla_x \phi(\tau, X_i(\tau)) \big) \diff \tau\\
&= m_i \sum_{i=1}^N \int_0^t V_i(\tau) \cdot \frac{\diff}{\diff \tau} \phi(\tau, X_i(\tau)) \diff \tau \\
&= \int_{\R^d}\phi(t, x)\cdot u(t,x)\diff\rho_t(x) - \int_{\R^d}\phi(0, x)\cdot u_0(x)\diff\rho_0(x) - m_i \sum_{i=1}^N \int_0^t \dot{V_i}(\tau) \cdot \phi(\tau, X_i(\tau))  \diff \tau.
\end{align*}
In order to obtain \eqref{eq:def_weak_moment}, it remains to notice that 
\begin{align*}
&-  m_i \sum_{i=1}^N \int_0^t \dot{V_i}(\tau) \cdot \phi(\tau, X_i(\tau))  \diff \tau =\\ &=\int_0^t \sum_{i=1}^N \sum_{j=1}^N  \int_0^tm_i\, m_j\, \nabla K_{\eps}(X_i(t)-X_j(t)) \cdot \phi(\tau, X_i(\tau))  \diff \tau + \int_0^t \sum_{i=1}^N m_i \, \kappa \,V_i(t)\cdot \phi(\tau, X_i(\tau))  \diff \tau\\
&= \int_0^t \int_{\R^d} \big( \nabla_x(K_\varepsilon\star\diff\rho_{\tau})\cdot\phi  +  \kappa\,u\cdot\phi \big) \diff\rho_\tau(x)\diff \tau.
\end{align*}
Finally, to see \eqref{ineq:def_weak_energy}, we multiply by $m_i\,V_i(t)$ the equation on $\dot{V_i}(t)$ in \eqref{eq:SPH_Rd} and we sum over $i=1, ..., N$  to obtain
$$
\frac{1}{2} \, \frac{\diff}{\diff \tau} \sum_{i=1}^N m_i \, |V_i|^2 = - \sum_{i,j=1}^N m_i\, m_j\, \nabla K_{\eps}(X_i(\tau)-X_j(\tau)) \cdot V_i(\tau)  - \kappa \sum_{i=1}^N m_i\, |V_i|^2
$$
By the symmetry of $K_{\eps}$, $\nabla K_{\eps}(X_i(\tau) - X_j(\tau)) = -\nabla K_{\eps}(X_j(\tau) - X_i(\tau))$ so that we recognize 
\begin{multline*}
 \sum_{i,j=1}^N m_i\, m_j\, \nabla K_{\eps}(X_i(\tau)-X_j(\tau)) \cdot V_i(\tau) = \frac{1}{2} \sum_{i,j=1}^N m_i\, m_j\, \nabla K_{\eps}(X_i(t) - X_j(t))\cdot (V_i(t) - V_j(t)) =\\
= \frac{1}{2} \frac{\diff}{\diff t} \int_{\R^d} K_{\eps} \star \rho_t \diff \rho_t = \frac{1}{2} \frac{\diff}{\diff t} \int_{\R^d} |L_{\eps} \star \rho_t|^2 \diff x.
\end{multline*}
Hence, we obtain \eqref{ineq:def_weak_energy} as an equality with $\mu = 0$. The proof of \eqref{eq:def_weak_cont_per}--\eqref{ineq:def_weak_energy_per} is analogous. 
\end{proof}

\section{Proof of Theorem \ref{thm:main}}\label{sect:proof_theorem_quadratic_pressure}

\noindent \underline{Step 1. Proof of \eqref{eq:energy_relative_main_bound_Rd}, \eqref{eq:energy_relative_main_bound_Rd_gronwall_estimate}, and \eqref{eq:energy_relative_main_bound_Rd_gronwall_estimate_negative_k}.} Note that since $(r, U)$ solves \eqref{sys:loc_euler}, it also solves
    \begin{equation}\label{sys:loc_euler_strong}
    \begin{aligned}
        \partial_tr + \DIV_x(r\, U) &= 0 &&\quad\mbox{in }(0,T)\times \R^d,\\
        \partial_t (r\,U)  + \DIV_x (r\, U \otimes U ) + r\nabla_x(K_\varepsilon\star r) + \kappa\,r U  &= e &&\quad\mbox{in }(0,T)\times \R^d, \\
    \end{aligned}
    \end{equation}
    where the error $e$ is defined as follows
    $$
        e := r\nabla_x(K_\varepsilon\star r) - r\nabla_x r.
    $$
    Moreover, for further use let use notice that since $\nabla^2_x r\in L^\infty((0, T)\times \R^d)$, we can estimate
    \begin{align}\label{ineq:error_estimate}
         \left|\frac{e}{r}\right| \leq \|\nabla_x^2 r\|_\infty \int_{\R^d}K_{\varepsilon}(x)|x|\diff x = \varepsilon \|\nabla^2_x r\|_\infty\int_{\R^d}K(x)|x|\diff x = \varepsilon\, C(m, d, \|\nabla^2_x r\|_\infty).
    \end{align}
    Moving forward we begin deriving the relative entropy inequality. First, we test \eqref{eq:def_weak_cont} by $\frac{1}{2}|U|^2$ and $K_\varepsilon\star r$ to obtain respectively
    \begin{equation}\label{eq:relative_U_2}
        \begin{split}
            \int_{\R^d}\frac{1}{2}|U|^2(t, x)\diff\rho^\varepsilon_t(x) - \int_{\R^d}\frac{1}{2}|U_0|^2\diff \rho^\varepsilon_0(x) = \int_0^t\int_{\R^d} \big(U\cdot \partial_t U +  u^\varepsilon\cdot(\nabla_x U)^{\top} U\big) \diff\rho^\varepsilon_\tau(x)\diff\tau
        \end{split}
    \end{equation}
    and
    \begin{equation}\label{eq:relative_K_r}
        \begin{split}
            \int_{\R^d}(K_\varepsilon\star r)(t, x)\diff\rho^\varepsilon _t(x) - \int_{\R^d}K_{\varepsilon}\star r_0\diff\rho^\varepsilon_0(x) = \int_0^t\int_{\R^d}\big(\partial_t(K_\varepsilon\star r) +  u^\varepsilon\cdot\nabla_xK_\varepsilon\star r\big) \diff\rho^\varepsilon_\tau(x)\diff\tau.
        \end{split}
    \end{equation}
    Furthermore, testing \eqref{eq:def_weak_moment} by $U$ gives us
    \begin{equation}\label{eq:relative_U}
        \begin{split}
            &\int_{\R^d} u^\varepsilon\cdot U\diff\rho^\varepsilon_t - \int_{\R^d}u^\varepsilon_0\cdot U_0\diff\rho^\varepsilon_0 = \int_0^t\int_{\R^d}\big(u^\varepsilon\cdot\partial_t U +  u^\varepsilon\otimes u^\varepsilon : \nabla_x U - \kappa\,u^\varepsilon\cdot U\big) \diff\rho_\tau^\varepsilon(x)\diff\tau\\
            &\qquad \qquad- \int_0^t\int_{\R^d}\nabla_x K_\varepsilon\star\diff\rho^\varepsilon_\tau\cdot U\diff \rho^\varepsilon_\tau(x)\diff\tau + \int_0^t\int_{\R^d}\nabla_x U:\diff\mu^\varepsilon_\tau(x)\diff \tau.
        \end{split}
    \end{equation}
    By adding together \eqref{eq:relative_U_2} and $\frac{1}{2}$ of the energy inequality \eqref{ineq:def_weak_energy}, and subtracting from the result \eqref{eq:relative_U} as well as \eqref{eq:relative_K_r}, we get (note that in the first line we included an additional term $\frac{1}{2}\int_{\R^d}|L_{m, \varepsilon}\star r_t|^2\diff x - \frac{1}{2}\int_{\R^d}|L_{m, \varepsilon}\star r_0|^2\diff x $ to form $\mathcal{E}^\varepsilon_{rel}(t)  - \mathcal{E}^\varepsilon_{rel}(0)$)
    \begin{equation}\label{ineq:relative_entropy_1}
        \begin{split}
            \mathcal{E}^\varepsilon_{rel}(t) + \frac{1}{2}\mathrm{tr}(\mu^\varepsilon_t(\R^d))  - \mathcal{E}^\varepsilon_{rel}(0) &\leq -\kappa \int_0^t\int_{\R^d}|u^\varepsilon|^2\diff\rho^\varepsilon_\tau(x)\diff\tau + \frac{1}{2}\int_{\R^d}|L_{m, \varepsilon}\star r_t|^2\diff x - \frac{1}{2}\int_{\R^d}|L_{m, \varepsilon}\star r_0|^2\diff x \\
            &\quad +\int_0^t\int_{\R^d} \big(U\cdot \partial_t U + u^\varepsilon\cdot(\nabla_x U)^{\top} U\big) \diff\rho^\varepsilon_\tau(x)\diff\tau\\
            &\quad - \int_0^t\int_{\R^d} \big(u^\varepsilon\cdot\partial_t U + u^\varepsilon\otimes u^\varepsilon : \nabla_x U - k\, u^\varepsilon\cdot U - \nabla_x K_\varepsilon\star\diff\rho^\varepsilon\cdot U\big) \diff \rho^\varepsilon_\tau(x)\diff\tau\\
            &\quad - \int_0^t\int_{\R^d}\nabla_x U : \diff\mu^\varepsilon_\tau(x)\diff\tau\\
            &\quad -\int_0^t\int_{\R^d}\big(\partial_t(K_\varepsilon\star r) +  u^\varepsilon\cdot\nabla_xK_\varepsilon\star r\big) \diff\rho^\varepsilon_\tau(x) \diff\tau.
        \end{split}
    \end{equation}
    Now, since $r > 0$, we deduce from \eqref{sys:loc_euler_strong} that the strong velocity $U$ satisfies the equation
    $$
    \partial_t U + (U\cdot\nabla_x )U + \kappa\,U + \nabla_x K_\varepsilon\star r = \frac{e}{r}.
    $$
 This allows to write some terms in \eqref{ineq:relative_entropy_1} as
    \begin{equation*}
        \begin{split}
            &\int_0^t\int_{\R^d} \big(U\cdot \partial_t U +  u^\varepsilon\cdot(\nabla_x U )^{\top} U\big) \diff\rho^\varepsilon_\tau(x)\diff\tau - \int_0^t\int_{\R^d} \big( u^\varepsilon\cdot\partial_t U + u^\varepsilon\otimes u^\varepsilon : \nabla_x U\big) \diff\rho^\varepsilon_\tau(x)\diff \tau \\
            & = \int_0^t\int_{\R^d} \big( (U - u^\varepsilon)\, \partial_t U + (U - u^\varepsilon)\cdot\nabla_x U\, U + (u^\varepsilon - U)\otimes (U - u^\varepsilon):\nabla_x U\big)\diff\rho^\varepsilon_\tau(x)\diff\tau\\
            & = \int_0^t\int_{\R^d}\Big(\kappa\, (u^\varepsilon - U)\cdot U + (u^\varepsilon - U)\cdot\nabla_x K_\varepsilon\star r - (u^\varepsilon - U)\cdot\frac{e}{r}\Big) \diff \rho^\varepsilon_\tau(x)\diff\tau\\
            &\qquad + \int_0^t\int_{\R^d}(U - u^\varepsilon)\otimes (u^\varepsilon - U):\nabla_x U\diff\rho^\varepsilon_\tau(x)\diff\tau,
        \end{split}
    \end{equation*}
    which combined with \eqref{ineq:relative_entropy_1} grants us
    \begin{equation}\label{ineq:relative_entropy_2}
        \begin{split}
            \mathcal{E}^\varepsilon_{rel}(t)  + \frac{1}{2}\mathrm{tr}(\mu^\varepsilon_t(\R^d))  - \mathcal{E}^\varepsilon_{rel}(0) &\leq -\kappa \int_0^t\int_{\R^d}|u^\varepsilon - U|^2\diff\rho^\varepsilon_\tau(x)\diff\tau + \frac{1}{2}\int_{\R^d}|L_{m, \varepsilon}\star r|^2\diff x - \frac{1}{2}\int_{\R^d}|L_{m, \varepsilon}\star r_0|^2\diff x \\
            &\quad -\int_0^t\int_{\R^d}\partial_t(K_\varepsilon\star r) +  u^\varepsilon\cdot\nabla_xK_\varepsilon\star r -  U\cdot\nabla_x K_\varepsilon\star\diff\rho^\varepsilon_\tau\diff\rho^\varepsilon_\tau(x)\diff\tau\\
            &\quad +\int_0^t\int_{\R^d}(u^\varepsilon - U)\cdot\nabla_x K_\varepsilon\star r + (u^\varepsilon - U)\otimes (U - u^\varepsilon) : \nabla_x U\diff\rho^\varepsilon_\tau(x)\diff\tau\\
            &\quad -\int_0^t\int_{\R^d}(u^\varepsilon - U)\cdot\frac{e}{r}\diff\rho^\varepsilon_\tau(x)\diff\tau  - \int_0^t\int_{\R^d}\nabla_x U : \diff\mu^\varepsilon_\tau(x)\diff\tau.
        \end{split}
    \end{equation}
    Now let us focus on the terms with the kernels $K_\varepsilon$ and $L_{m, \varepsilon}$. Observe that
    \begin{equation}\label{eq:chain_rule_in_time_derving_inequality_Rd}
        \frac{1}{2}\int_{\R^d}|L_{m, \varepsilon}\star r|^2\diff x - \frac{1}{2}\int_{\R^d}|L_{m, \varepsilon}\star r_0|^2\diff x = \int_0^t\int_{\R^d}r\, \partial_t(K_\varepsilon\star r)\diff x\diff\tau.
    \end{equation}
    Thus, adding and subtracting the same expressions to \eqref{eq:chain_rule_in_time_derving_inequality_Rd}, we deduce that some terms in \eqref{ineq:relative_entropy_2} can be written as follows:
    \begin{equation}\label{eq:relative_reduce_K_r}
        \begin{split}
            &\frac{1}{2}\int_{\R^d}|L_{m, \varepsilon}\star r|^2\diff x - \frac{1}{2}\int_{\R^d}|L_{m, \varepsilon}\star r_0|^2\diff x -\int_0^t\int_{\R^d} \big(\partial_t(K_\varepsilon\star r) + u^\varepsilon\cdot\nabla_xK_\varepsilon\star r -  U\cdot\nabla_x K_\varepsilon\star\diff\rho^\varepsilon_\tau \big) \diff\rho^\varepsilon_\tau(x)\diff\tau\\
            &\quad +\int_0^t\int_{\R^d}(u^\varepsilon - U)\cdot\nabla_x K_\varepsilon\star r\diff\rho^\varepsilon_\tau(x)\diff\tau = -\int_0^t\int_{\R^d} \partial_t(K_\varepsilon\star r)\diff(\rho^\varepsilon_\tau - r_\tau)(x)\diff\tau\\
            &\qquad+\int_0^t\int_{\R^d}U\cdot \nabla_xK_{\varepsilon}\star(\rho^\varepsilon_\tau - r)\diff\rho^\varepsilon_\tau(x)\diff\tau.
        \end{split}
    \end{equation}
    Using the continuity equation for the strong density $r$, and integrating by parts we have
    \begin{align*}
        &-\int_0^t\int_{\R^d}\partial_t(K_\varepsilon\star r)\diff(\rho^\varepsilon_\tau - r_\tau)(x)\diff\tau = -\int_0^t\int_{\R^d}K_\varepsilon\star(\diff\rho^\varepsilon_\tau - \diff r_\tau)\partial_t r\diff x\diff\tau\\
        &\qquad\qquad= \int_0^t\int_{\R^d}K_\varepsilon\star(\diff\rho^\varepsilon_\tau - \diff r_\tau)\DIV_x(rU)\diff x\diff\tau = -\int_0^t\int_{\R^d}r U\cdot\nabla_x K_\varepsilon\star(\diff \rho^\varepsilon_\tau - \diff r_\tau)\diff x\diff\tau.
    \end{align*}
    Combined with \eqref{eq:relative_reduce_K_r} this reduces \eqref{ineq:relative_entropy_2} to
    \begin{equation}\label{ineq:relative_entropy_3}
        \begin{split}
            \mathcal{E}^\varepsilon_{rel}(t) + \frac{1}{2}\mathrm{tr}(\mu^\varepsilon_t(\R^d))  - \mathcal{E}^\varepsilon_{rel}(0) &\leq -\kappa \int_0^t\int_{\R^d}|u^\varepsilon - U|^2\diff\rho^\varepsilon_\tau(x)\diff\tau\\
            &\quad+ \int_0^t\int_{\R^d}(u^\varepsilon - U)\otimes (U - u^\varepsilon) : \nabla_x U\diff\rho^\varepsilon_\tau(x)\diff\tau\\
            &\quad + \int_0^t\int_{\R^d}U\cdot\nabla_xK_\varepsilon\star(\diff \rho^\varepsilon_\tau - \diff r_\tau)\diff(\rho^\varepsilon_\tau -  r_\tau)(x)\diff\tau\\
            &\quad -\int_0^t\int_{\R^d}(u^\varepsilon - U)\cdot\frac{e}{r}\diff\rho^\varepsilon_\tau(x)\diff\tau\\
            &\quad - \int_0^t\int_{\R^d}\nabla_x U : \diff\mu^\varepsilon_\tau(x)\diff\tau.
        \end{split}
    \end{equation}
    Having this, we may start obtaining proper bounds on the right-hand side in term of the relative entropy $\mathcal{E}^\varepsilon_{rel}$. For the second term we estimate
    \begin{align}\label{ineq:relative_bound_1}
    \left|\int_0^t\int_{\R^d}(u^\varepsilon - U)\otimes (U - u^\varepsilon) : \nabla_x U\diff\rho^\varepsilon_\tau(x)\diff\tau\right|\leq \|\nabla_x U\|_\infty\int_0^t\int_{\R^d}|u^\varepsilon - U|^2\diff\rho^\varepsilon_\tau(x)\diff\tau.
    \end{align}
    For the fourth term we apply the inequality on the error term \eqref{ineq:error_estimate}, the fact that $\rho^\eps$ is a probability measure and Cauchy-Schwarz inequality to see that
    \begin{equation}\label{ineq:relative_bound_4}
\begin{split}
        \left|\int_0^t\int_{\R^d}(u^\varepsilon - U)\cdot\frac{e}{r}\diff\rho^\varepsilon_\tau(x)\diff\tau\right|&\leq \frac{1}{2}\int_0^t\int_{\R^d}|u^\varepsilon - U|^2\diff\rho^\varepsilon_\tau(x)\diff\tau + \frac{1}{2}\int_0^t\int_{\R^d}\left|\frac{e}{r}\right|^2\diff\rho^\varepsilon_\tau(x)\diff\tau\\
        &\leq \frac{1}{2} \int_0^t  \mathcal{E}^\varepsilon_{rel}(\tau)\diff \tau +  \varepsilon^2\,  C(T, m, d, \|\nabla^2_x r\|_\infty).
    \end{split}
\end{equation}
    Furthermore, for the term involving the concentration measure $\mu$ we can bound using Proposition \ref{prop:bound_symmetric_measure_trace}
    \begin{align}\label{ineq:relative_bound_5}
        \left|\int_0^t\int_{\R^d}\nabla_x U : \diff\mu^\varepsilon_\tau(x)\diff\tau\right| \leq C\|\nabla_x U\|_\infty\int_0^t\mathrm{tr}(\mu^\varepsilon_\tau(\R^d))\diff\tau.
    \end{align}
    Hence, let us focus on the term concerning the kernel $K_\varepsilon$. We rewrite
    \begin{equation}\label{eq:relative_bound_2}
        \begin{split}
            &\int_0^t\int_{\R^d}U\cdot\nabla_xK_\varepsilon\star(\diff\rho^\varepsilon_\tau - \diff r_\tau)\diff(\rho^\varepsilon_\tau - \diff r_\tau)(x)\diff\tau\\
            &= \int_0^t\int_{\R^d}L_{m, \varepsilon}\star(U(\diff\rho^\varepsilon_\tau - \diff r_\tau))\cdot\nabla_x L_{m, \varepsilon}\star(\diff\rho^\varepsilon_\tau - \diff r_\tau)\diff x\diff\tau \\
            & = \int_{0}^t\int_{\R^d}\left(L_{m, \varepsilon}\star(U(\diff\rho^\varepsilon_\tau - \diff r_\tau)) - U\,L_{m, \varepsilon}\star(\diff\rho^\varepsilon
             - \diff r)\right)\cdot\nabla_x L_{m, \varepsilon}\star(\diff\rho^\varepsilon_\tau - \diff r_\tau)\diff x\diff\tau\\
             &\quad + \int_0^t\int_{\R^d}L_{m, \varepsilon}\star(\diff\rho^\varepsilon
             - \diff r)U\cdot\nabla_x L_{m, \varepsilon}\star(\diff\rho^\varepsilon_\tau - \diff r_\tau)\diff x\diff\tau\\
              & = \int_{0}^t\int_{\R^d}\left(L_{m, \varepsilon}\star(U(\diff\rho^\varepsilon_\tau - \diff r_\tau)) - U\,L_{m, \varepsilon}\star(\diff\rho^\varepsilon
             - \diff r)\right)\cdot\nabla_x L_{m, \varepsilon}\star(\diff\rho^\varepsilon_\tau - \diff r_\tau)\diff x\diff\tau\\
             &\quad + \frac{1}{2}\int_0^t\int_{\R^d}U\cdot\nabla_x |L_{m, \varepsilon}\star(\diff\rho^\varepsilon_\tau - \diff r_\tau)|^2\diff x\diff\tau\\
             & = \int_{0}^t\int_{\R^d}\left(L_{m, \varepsilon}\star(U(\diff\rho^\varepsilon_\tau - \diff r_\tau)) - U\,L_{m, \varepsilon}\star(\diff\rho^\varepsilon
             - \diff r)\right)\cdot\nabla_x L_{m, \varepsilon}\star(\diff\rho^\varepsilon_\tau - \diff r_\tau)\diff x\diff\tau\\
             &\quad - \frac{1}{2}\int_0^t\int_{\R^d}\DIV_xU\,|L_{m, \varepsilon}\star(\diff\rho^\varepsilon_\tau - \diff r_\tau)|^2\diff x\diff\tau\\
             & = \int_0^t\int_{\R^d}\int_{\R^d}\nabla L_{m, \varepsilon}(x-y)\cdot (U(x) - U(y))L_{m, \varepsilon}\star(\diff\rho^\varepsilon_\tau - \diff r_\tau)(x)\diff(\rho^\varepsilon_\tau - r_\tau)(y)\diff x\diff\tau \\
             &\quad+ \int_0^t\int_{\R^d}\DIV_xU|L_{m, \varepsilon}\star(\diff\rho^\varepsilon_\tau - \diff r_\tau)|^2\diff x\diff\tau - \frac{1}{2}\int_0^t\int_{\R^d}\DIV_xU\,|L_{m, \varepsilon}\star(\diff\rho^\varepsilon_\tau - \diff r_\tau)|^2\diff x\diff\tau\\
             &= \int_0^t\int_{\R^d}\int_{\R^d}\nabla L_{m, \varepsilon}(x-y)\cdot (U(x) - U(y))L_{m, \varepsilon}\star(\diff\rho^\varepsilon_\tau - \diff r_\tau)(x)\diff(\rho^\varepsilon_\tau - r_\tau)(y)\diff x\diff\tau \\
             &\quad+ \frac{1}{2}\int_0^t\int_{\R^d}\DIV_xU\,|L_{m, \varepsilon}\star(\diff\rho^\varepsilon_\tau - \diff r_\tau)|^2\diff x\diff\tau.
        \end{split}
    \end{equation}
Since the second term can be easily treated, we focus on the first one. Using Taylor's formula with Lagrange's remainder we can write
    \begin{align*}
        U_i(x) - U_i(y) = \nabla_x U_i(x)\cdot (x - y) + \frac{1}{2}\nabla^{2}U_i(\zeta_{x, y}):(x - y)\otimes (x-y),
    \end{align*}
    Hence, we may estimate using H\"{o}lder's inequality
    \begin{equation}\label{ineq:taylor_formula_lemma}
        \begin{split}
            &\left|\int_{\R^d}\int_{\R^d}\nabla L_{m, \varepsilon}(x-y)\cdot (U(x) - U(y))L_{m, \varepsilon}\star(\diff\rho^\varepsilon_\tau - \diff r_\tau)(x)\diff(\rho^\varepsilon_\tau - r_\tau)(y)\diff x\right|\\
            &\quad\leq \|\nabla_x U\cdot (\diff\rho^\varepsilon_\tau - \diff r_\tau)\star(\nabla L_{m, \varepsilon}(x)\,x)\|_2\|L_{m, \varepsilon}\star(\diff\rho^\varepsilon_\tau - \diff r_\tau)\|_2 \\
            &\qquad+ \frac{1}{2}\int_{\R^d}\int_{\R^d}|\nabla L_{m, \varepsilon}(x-y)||x-y|^2 |\nabla^{2}U(\zeta_{x, y})|L_{m, \varepsilon}\star(\diff\rho^\varepsilon_\tau - \diff r_\tau)(x)\diff |\rho^\varepsilon_\tau - r_\tau|(y)\diff x\\
            &\quad\leq \|\nabla_x U\|_\infty \|(\diff\rho^\varepsilon_\tau - \diff r_\tau)\star(\nabla L_{m, \varepsilon}(x)\,x)\|_2\|L_{m, \varepsilon}\star(\diff\rho^\varepsilon_\tau - \diff r_\tau)\|_2 \\
            &\qquad+ \frac{1}{2}\int_{\R^d}\int_{\R^d}|\nabla L_{m, \varepsilon}(x-y)||x-y|^2 |\nabla^{2}U(\zeta_{x, y})||L_{m, \varepsilon}\star(\diff\rho^\varepsilon_\tau - \diff r_\tau)|(x)\diff |\rho^\varepsilon_\tau - r_\tau|(y)\diff x\\
            &\quad =: I_1 + I_2.
        \end{split}
    \end{equation}
    Let us handle the $I_1$ first. By the isometry of Fourier transform in $L^2(\R^d)$ and its standard properties, we can deduce for an arbitrary, finite measure $\nu\in \mathcal{M}(\R^d)$
    \begin{equation}\label{ineq:l2_bound_conv_lemma}
    \begin{split}
        &\|(\nabla L_{m, \varepsilon}(x)\,x)\star\nu\|^2_2 = \sum_{i = 1}^d\|(\partial_{x_i} L_{m, \varepsilon}(x)\,x) \star\nu\|^2_2 = \sum_{i = 1}^d \int_{\R^d}|\hat{\nu}(\xi)|^2\left| ({\partial_{x_i} L_{m, \varepsilon}(x) \, x})^{\wedge} (\xi)\right|^2\diff \xi \\
        &\lesssim \sum_{i = 1}^d\int_{\R^d}|\hat{\nu}(\xi)|^2\left|\nabla_{\xi} \left(\frac{\xi_i}{(1 + 4\pi^2|\varepsilon\xi|^2)^{m/2}}\right)\right|^2\diff \xi\\
        &\quad = \sum_{i = 1}^d\sum_{j = 1, j\neq i}^d\int_{\R^d}|\hat{\nu}(\xi)|^2\left|\frac{4\pi^2m\eps^2\xi_j\xi_i}{(1 + 4\pi^2|\eps\xi|^2)^m}\right|^2\diff \xi\\
&\quad \phantom{ = }  + \sum_{i = 1}^d\int_{\R^d}|\hat{\nu}(\xi)|^2\left|\frac{1}{(1 + 4\pi^2|\eps\xi|^2)^{m/2}} + \frac{4\pi^2m\eps^2|\xi_i|^2}{(1 + 4\pi^2|\eps\xi|^2)^m}\right|^2\diff\xi.
    \end{split}
    \end{equation}
    It is quite simple to verify that as long as $m\geq 2$
    \begin{align*}
        \left|\frac{4\pi^2m\eps^2\xi_j\xi_i}{(1 + 4\pi^2|\eps\xi|^2)^m}\right| \leq \frac{m}{(1 + 4\pi^2|\varepsilon \xi|^2)^{m/2}},\quad \left|\frac{1}{(1 + 4\pi^2|\eps\xi|^2)^{m/2}} + \frac{4\pi^2m\eps^2|\xi_i|^2}{(1 + 4\pi^2|\eps\xi|^2)^m}\right| \leq \frac{1+m}{(1 + 4\pi^2|\varepsilon \xi|^2)^{m/2}},
    \end{align*}
    which combined with \eqref{ineq:l2_bound_conv_lemma} grants us
    \begin{align*}
        \|\nabla L_{m, \varepsilon}(x)x\star\nu\|^2_2 \leq C(m, d) \int_{\R^d}|\hat{\nu}(\xi)|^2\left|\frac{1}{(1 + 4\pi^2|\varepsilon \xi|^2)^{m/2}}\right|^2\diff \xi =  C(m, d) \|L_{m, \varepsilon}\star\nu\|_2^2.
    \end{align*}
    We conclude
    \begin{align}\label{ineq:I_1_whole_domain}
         |I_1| \leq C(m, d)\|U\|_{W^{1, \infty}(\R^d)}\|L_{m, \varepsilon}\star(\diff\rho^\eps_\tau - \diff r_\tau)\|^2_2.
    \end{align}
    For $I_2$ we utilize Lemma \ref{lem:bound_for_nabla_bessel_potential} (to obtain that ${|\nabla L_{m, \varepsilon}(x)||x|^2} \leq C\, \varepsilon\, L_{m+4,\eps}(x)$), H\"{o}lder's inequality, Young's inequality and the energy bounds \eqref{ineq:def_weak_energy} to estimate
    \begin{equation}\label{ineq:I_2_whole_domain}
        \begin{split}
            &\int_{\R^d}\int_{\R^d}|\nabla L_{m, \varepsilon}(x-y)||x-y|^2 |\nabla^{2}U(\zeta_{x, y})|L_{m, \varepsilon}\star(\diff\rho^\varepsilon_\tau - \diff r_\tau)(x)\diff |\rho^\varepsilon_\tau - r_\tau|(y)\diff x\\
            &\quad \leq C\eps\|U\|_{W^{2,\infty}(\R^d)}\int_{\R^d}L_{m+4,\eps}\star(\diff\rho^\eps_\tau + \diff r_\tau)(x) |L_{m,\eps}\star(\diff\rho^\eps_\tau - \diff r_\tau)|(x)\diff x\\
            &\quad \leq C\eps\|U\|_{W^{2,\infty}(\R^d)}\|L_{4,\eps}\star L_{m,\eps}\star(\diff\rho^\eps_\tau + \diff r_\tau)\|_2 \|L_{m,\eps}\star(\diff\rho^\eps_\tau - \diff r_\tau)\|_2\\
            &\quad\leq C\eps\|U\|_{W^{2,\infty}(\R^d)}(\|L_{m,\eps}\star \diff\rho^\eps_\tau\|_2^2 +\|L_{{m,\eps}}\star r_\tau\|_2^2)\\
            &\quad\leq C(\|U\|_{L^\infty(0, T; W^{2,\infty}(\R^d))}, \|r\|_{L^\infty((0, T)\times \R^d}, \|u_0^\eps\|_{L^2(\R^d; \diff\rho_0^\eps)})\, \eps.
        \end{split}
    \end{equation}
    Combining \eqref{ineq:I_1_whole_domain} and \eqref{ineq:I_2_whole_domain} we get
    \begin{equation}\label{ineq:relative_bound_3}
        \begin{split}
            &\left|\int_0^t\int_{\R^d}\int_{\R^d}\nabla L_{m, \varepsilon}(x-y)\cdot (U(x) - U(y))L_{m, \varepsilon}\star(\diff\rho^\varepsilon_\tau - \diff r_\tau)(x)\diff (\diff\rho^\varepsilon_\tau - \diff r_\tau)(y)\diff x\diff\tau\right|\\
            &\leq C(\|U\|_{L^\infty(0, T; W^{1, \infty}(\R^d))}, m, d)\int_0^t\|L_{m, \varepsilon}\star(\diff\rho^\varepsilon_\tau
             - \diff r_\tau)\|_2^2\diff\tau\\
             &\qquad+ \varepsilon\, C(\|U\|_{L^\infty(0, T; W^{2,\infty}(\R^d))}, \|r\|_{L^\infty((0, T)\times \R^d}, \|u_0^\eps\|_{L^2(\R^d; \diff\rho_0^\eps)})\, T.
        \end{split}
    \end{equation}
    Combining \eqref{ineq:relative_entropy_3} with the obtained bounds \eqref{ineq:relative_bound_1}, \eqref{eq:relative_bound_2}, \eqref{ineq:relative_bound_3}, \eqref{ineq:relative_bound_4}, \eqref{ineq:relative_bound_5} grants us
    \begin{equation*}
        \begin{split}
            \mathcal{E}^\varepsilon_{rel}(t) + \frac{1}{2}\mathrm{tr}\left(\mu^\eps_t(\R^d)\right) &+ \kappa \int_0^t\int_{\R^d}|u^\varepsilon - U|^2\diff \rho^\varepsilon_\tau(x)\diff\tau \\
            &\leq \mathcal{E}^\varepsilon_{rel}(0) + C(\|U\|_{L^\infty(0, T; W^{1, \infty}(\R^d))}, m, d)\int_0^t  \mathcal{E}^\varepsilon_{rel}(\tau)+ \frac{1}{2}\mathrm{tr}\left(\mu^\eps_\tau(\R^d)\right) \diff \tau \\
             &\qquad + \varepsilon \, C(\|U\|_{L^\infty(0, T; W^{2,\infty}(\R^d))}, \|r\|_{L^\infty((0, T)\times \R^d}, \|u_0^\eps\|_{L^2(\R^d; \diff\rho_0^\eps)})\, T\\
             &\qquad + \varepsilon^2 \, C(K, T, m, d, \|\nabla^2_x r\|_\infty).
        \end{split}
    \end{equation*}
When $\kappa \geq 0$, we can apply Gr\"{o}nwall's inequality to obtain \eqref{eq:energy_relative_main_bound_Rd_gronwall_estimate}. On the other hand, if $\kappa < 0$, we estimate it by $ -\kappa \int_0^t\int_{\R^d}|u^\varepsilon - U|^2\diff \rho^\varepsilon_\tau(x)\diff\tau \leq -2\,\kappa \int_0^t \mathcal{E}^\varepsilon_{rel}(\tau)\diff \tau$. Applying Gr\"{o}nwall's inequality yields \eqref{eq:energy_relative_main_bound_Rd_gronwall_estimate_negative_k}.\\
 
\noindent \underline{Step 2. Proof of \eqref{eq:final_estimate_W_2_on_R_d}, \eqref{eq:final_estimate_flat_norm_on_R_d}, \eqref{eq:final_estimate_neg_SS_on_R_d}, and of the convergence statements.} 
We note that 
    \begin{align}\label{conv:from_relative_entropy}
        \int_{\R^d}|u^\varepsilon - U|^2(t)\diff \rho^\varepsilon_t(x) + \int_{\R^d}|L_{m, \varepsilon}\star(\diff\rho^\varepsilon - \diff r)|^2(t)\diff x + \mathrm{tr}(\mu^\varepsilon_t(\R^d)) = 2\, \mathcal{E}_{rel}^\varepsilon(t)+ \mathrm{tr}(\mu^\varepsilon_t(\R^d))
    \end{align}
    for a.e. $t\in (0, T)$. Therefore, \eqref{eq:final_estimate_W_2_on_R_d} and  \eqref{eq:final_estimate_flat_norm_on_R_d} follow from Lemma \ref{lem:estimate_wasserstein_flat} and \eqref{eq:energy_relative_main_bound_Rd_gronwall_estimate} (or \eqref{eq:energy_relative_main_bound_Rd_gronwall_estimate_negative_k} if $\kappa < 0$).  Regarding \eqref{eq:final_estimate_neg_SS_on_R_d}, by Proposition~\ref{prop:bound_symmetric_measure_trace}, we have
    $
    \|\mu^\varepsilon_t\|_{TV}\leq C\, \mathrm{tr}(\mu_t^\varepsilon(\R^d)).
    $
    Similarly, the estimate in $H^{-m}(\R^d)$ norm comes from the energy $\mathcal{E}_{rel}^\varepsilon(t)$, as for any $F\in H^m(\R^d)$
    \begin{equation}\label{eq:nonlocal_operator_gives_convergence_in_neg_SS_linear_pressure}
        \begin{split}
            &\left|\int_{\R^d}F(x)\diff(\rho^\varepsilon_t - r_t)(x)\right|^2 = \left|\int_{\R^d}\hat{F}(-\xi)\widehat{(\rho^\varepsilon - r)}(\xi)\diff\xi\right|^2\\
            &\leq \|(1+4\pi^2|\xi|^2)^{m/2}|\hat{F}(\xi)|^2\|_2^2\int_{\R^d}\frac{|\widehat{(\rho^\varepsilon - r)}(\xi)|^2}{(1 + 4\pi^2|\xi|^2)^{m}}\diff \xi \leq C\|F\|_{H^m}^2 \int_{\R^d}|L_{m, \varepsilon}\star(\diff\rho^\varepsilon_t - \diff r_t)|^2\diff x.
        \end{split}
    \end{equation}
    Hence, for a numerical constant $C$, $
        \|\rho^\varepsilon - r\|_{H^{-m}(\R^d)} + \|\mu^\varepsilon_t\|_{TV} \leq C\, (\mathcal{E}_{rel}^\varepsilon(t)+ \mathrm{tr}(\mu^\varepsilon_t(\R^d)))$ and we obtain \eqref{eq:final_estimate_neg_SS_on_R_d}.  The convergence statements follow now directly from \eqref{eq:energy_relative_main_bound_Rd_gronwall_estimate} (or \eqref{eq:energy_relative_main_bound_Rd_gronwall_estimate_negative_k} if $\kappa < 0$), \eqref{eq:final_estimate_W_2_on_R_d}, \eqref{eq:final_estimate_flat_norm_on_R_d}, and \eqref{eq:final_estimate_neg_SS_on_R_d}.

\section{Proof of Theorem \ref{thm:main_per}}\label{sect:proof_theorem_general_pressure} 

\subsection{An auxilliary decomposition}\label{sect:auxilliary_results_periodic_case} We first consider the quantity which will appear in the estimates:
    \begin{align*}
        \mathcal{A}(x,y) := \nabla_x L_{m, \varepsilon}^{per}(x-y)\cdot (U(x) - U(y)) + \DIV_x U(x)\, L_{m,\eps}^{per}(x-y)
    \end{align*}
and the related convolution operator $\int_{\T^d} \mathcal{A}(x,y) \diff \nu(y)$.
\begin{lem}\label{lem:operator-mathcal-A-decomposition}
Suppose that $\nu\in \mathcal{M}(\T^d)$, $1 < p < +\infty$ and $m > d$.
Then, there exist kernels $P_{i,j}$ and a function $S_2$ such that 
\begin{equation}\label{eq:final_decomposition_mathcal_A}
\int_{\T^d} \mathcal{A}(x,y) \diff \nu(y) = \sum_{j = 1}^d\partial_{x_j}U_i(x)\sum_{i = 1}^d P_{ij} \ast \nu (x) + \int S_2(x,y) \diff \nu(y).
\end{equation}
Moreover, for all $p \in (1,\infty)$ and $m>d$, there exists a constant $C(p, d, m) > 0$, such that for any $i, j = 1, ..., d$
\begin{equation}\label{eq:bound_P_ij}
    \|P_{ij}\star \nu\|_p \leq C\|L_{m, \varepsilon}^{per}\star \nu\|_p.
\end{equation}
Furthemore, there exists a constant $C(m,d)> 0$ such that for any $x,y\in\T^d$
\begin{equation}\label{eq:bound_on_the_remainder_S2}
        |S_2(x,y)|\leq C\eps\|U\|_{W^{2,\infty}(\T^d)} L_{m+4, \eps}^{per}(x- y).
\end{equation}
\end{lem}
\begin{proof}
\underline{Step 1. Decomposition of $\mathcal{A}$.} We will prove
\begin{equation}\label{eq:first_splitting_mathcal_A}
    \begin{split}
        \mathcal{A}(x, y)&= \sum_{i\neq j}^d\partial_{x_j}U_i(x)\sum_{l\in\I^d}\partial_{x_i} L_{m, \varepsilon}(x+l-y)(x+ l - y)_j\\
        &\qquad+ \sum_{i = 1}^d\sum_{l\in\I^d}\partial_{x_i}U_i(x)(\partial_{x_i}L_{m, \varepsilon}(x+l-y)(x+l-y)_i + L_{m,\eps}(x+l-y))\\
        &\qquad + \sum_{i = 1}^d\sum_{l\in\I^d}\partial_{x_i} L_{m, \varepsilon}(x+l-y)\nabla^2_x U_i(\zeta_{x, y, l}) : (x+ l- y)\otimes (x+l -y) =: S^1_1 + S_1^2 + S_2.
    \end{split}
\end{equation}
Using the definition of the periodization of $L_{m, \varepsilon}$ we get
\begin{equation*}
    \begin{split}
        \nabla_x L_{m, \varepsilon}^{per}(x-y)\cdot (U(x) - U(y)) = \sum_{i=1}^d\partial_{x_i} L_{m, \varepsilon}^{per}(x-y)(U_i(x) - U_i(y)) =\sum_{i = 1}^d \sum_{l\in \I^d}\partial_{x_i} L_{m, \varepsilon}(x+l-y)(U_i(x) - U_i(y)).
    \end{split}
\end{equation*}
Since $U$ is $l$-periodic for any $l\in \I^d$, we can identify $U_i(x) = U_i(x+l)$ and utilize the Taylor expansion to obtain
\begin{equation*}
    \begin{split}
        &\sum_{i = 1}^d\sum_{l\in\I^d}\partial_{x_i} L_{m, \varepsilon}(x+l-y)(U_i(x) - U_i(y)) = \sum_{i = 1}^d\sum_{l\in\I^d}\partial_{x_i} L_{m, \varepsilon}(x+l-y)(U_i(x+l) - U_i(y))\\
        &= \sum_{i = 1}^d\sum_{l\in\I^d}\partial_{x_i} L_{m, \varepsilon}(x+l-y)\nabla_x U_i(x+l)\cdot (x+ l - y) \\
        &\qquad + \sum_{i = 1}^d\sum_{l\in\I^d}\partial_{x_i} L_{m, \varepsilon}(x+l-y)\nabla^2_x U_i(\zeta_{x,y, l}) : (x+ l - y)\otimes (x+l -y).
    \end{split}
\end{equation*}
Noting that $\nabla_x U_i(x+l) = \nabla_x U_i$ and expanding the scalar product, we easily obtain the claim.\\
\noindent \underline{Step 2. The first part of the decomposition ($S^1_1 + S_1^2$) and the bound \eqref{eq:bound_P_ij}. } We will show that with $S^1_1$, $S^2_1$ as in \eqref{eq:first_splitting_mathcal_A} we have $S^1_1 + S_1^2 = \sum_{j = 1}^d\partial_{x_j}U_i(x)\sum_{i = 1}^d P_{ij} \ast \nu (x)$. Indeed, it is sufficient to define
\begin{align*}
P_{ij}(x) &:= \sum_{l\in\I^d}\partial_{x_i} L_{m, \varepsilon}(x+l)(x+ l)_j,  \qquad \qquad \qquad (i \neq j) \\
P_{ii}(x) &= \sum_{l\in\I^d}\partial_{x_i} L_{m, \varepsilon}(x+l)(x+ l)_i + L_{m, \varepsilon}(x+l).
\end{align*}
To prove the bound on the operators $P_{ij}$ we will use Fourier analysis. We observe that it is a periodization of the function $\partial_{x_i}L_{m, \varepsilon}\, x_j$ (or $\partial_{x_i}L_{m, \varepsilon}\, x_i + L_{m, \varepsilon}$ if $i=j$). Hence, by the Poisson's summation formula \cite[Chapter 1, \textsection\,10.3]{MR1800316} and the uniqueness of Fourier's coefficients, we have
$$
\hat{P}_{ij}(k) = \widehat{\partial_{x_i}L_{m, \varepsilon} x_j}(k) = -\partial_{\xi_j}(\xi_i\hat{L}_{m,\eps})(k) = -m\frac{4\pi^2\eps^2 k_i k_j}{(1 + 4\pi^2|\eps k|^2)^{\frac{m}{2} + 1}},\, (\mbox{if } i\neq j),
$$
$$
\hat{P}_{ii}(k) = \widehat{\partial_{x_i}L_{m, \varepsilon} x_i}(k)  + \widehat{L_{m, \varepsilon}}(k) = -\partial_{\xi_i}(\xi_i\hat{L}_{m,\eps})(k) + \hat{L}_{m,\eps}(k) = - m\frac{4\pi^2\eps^2 |k_i|^2}{(1 + 4\pi^2|\eps k|^2)^{\frac{m}{2} + 1}}.
$$
For $i\neq j$ we can rewrite equivalently
\begin{equation}\label{eq:computation_FT_of_Pij}
    -m\frac{4\pi^2\eps^2 k_i k_j}{(1 + 4\pi^2|\eps k|^2)^{\frac{m}{2} + 1}} = -m\frac{k_i}{|k|}\frac{k_j}{|k|}\left(\frac{1}{(1 + 4\pi^2|\eps k|^2)^{\frac{m}{2}}} -\frac{1}{(1 + 4\pi^2|\eps k|^2)^{\frac{m}{2} + 1}}\right).
\end{equation}
We recognize the coefficients in front of the expression as multipliers defining the Riesz transform
$$
R_i^{per}(f) := \sum_{k\in \I^d}m_i(k)\hat{f}(k)\exp(2\pi i k\cdot x),\, f\in L^p(\T^d),\, 1 < p < +\infty.
$$
with
$$
m_i(k) = -i\frac{k_i}{|k|},\, k\neq 0 \text{ and }m_i(0) = 0.
$$
On the other hand, the terms in the brackets in \eqref{eq:computation_FT_of_Pij} are $L^{per}_{2,\eps}(k)\, L_{m, \varepsilon}^{per}(k) - L_{m, \varepsilon}^{per}(k)$. Therefore, we can identify
\begin{equation}\label{eq:definition_pij_different_ij}
    P_{ij}\star \nu =-m R_i^{per}R_j^{per}(L^{per}_{2,\eps}\star L_{m, \varepsilon}^{per}\star \nu - L_{m, \varepsilon}^{per}\star \nu),\text{ for }\hat{L}_2(\xi) := \frac{1}{1 + 4\pi^2|\xi|^2},\, i\neq j,
\end{equation}
Similarly, we have
\begin{equation}\label{eq:definition_pij_same_ij}
    P_{ii}\star \nu =  -mR_i^{per}R_i^{per}((L^{per}_{2,\eps}\star L_{m, \varepsilon}^{per}\star \nu) - L_{m, \varepsilon}^{per}\star \nu).
\end{equation}
We are now ready to prove the estimate. Note that as $m > d$ we have $L_m\in L^1(\R^d)\cap L^\infty(\R^d)$, which means that all the terms are well defined. Hence, as $R_i^{per}$ is a bounded operator from $L^p(\T^d)$ to $L^p(\T^d)$ we can estimate using Young's inequality
    \begin{align*}
        \|P_{ij}\star \nu\|_p &\leq m\|R_i^{per}R_i^{per}((L^{per}_{2,\eps}\star L_{m, \varepsilon}^{per})\star \nu - L_{m, \varepsilon}^{per}\star \nu)\|_p\\
        &\leq Cm\|(L^{per}_{2,\eps}\star L_{m, \varepsilon}^{per})\star \nu\|_p + Cm \|L_{m, \varepsilon}^{per}\star \nu\|_p\\
        &\leq Cm\|L_{m, \varepsilon}^{per}\star \nu\|_p + Cm\|L_{2,\eps}^{per}\|_1\|L_{m, \varepsilon}^{per}\star \nu\|_p \leq 2Cm\|L_{m, \varepsilon}^{per}\star \nu\|_p.
    \end{align*}
\underline{Step 3. Proof of the estimate \eqref{eq:bound_on_the_remainder_S2} on $S_2$.} By Lemma \ref{lem:bound_for_nabla_bessel_potential}
    \begin{equation*}
        \begin{split}
            |S_2| & = \sum_{i = 1}^d\sum_{l\in\I^d}\frac{1}{\eps^{d-1}}\partial_{x_i} L_{m}\left(\frac{x+l-y}{\eps}\right)\nabla^2_x U_i(\zeta_{x, y, l}) : \left(\frac{x+ l- y}{\eps}\right)\otimes \left(\frac{x+ l- y}{\eps}\right)\\
            &\leq \|U\|_{W^{2,\infty}(\T^d)}\sum_{i=1}^d \sum_{l\in\I^d}\frac{1}{\eps^{d-1}}\left|\partial_{x_i}L_{m}\left(\frac{x+l - y}{\eps}\right)\right|\left|\frac{x+l - y}{\eps}\right|^2\\
            &\leq C\, \|U\|_{W^{2,\infty}(\T^d)}\sum_{l\in\I^d}\frac{1}{\eps^{d-1}}L_{m+4}\left(\frac{x+l-y}{\eps}\right) = C\eps \|U\|_{W^{2,\infty}(\T^d)}L^{per}_{m+4,\eps}(x-y).
        \end{split}
    \end{equation*}
\end{proof}
\subsection{Proof of Theorem \ref{thm:main_per}.} The proof is divided into several steps.\\

\noindent \underline{Step 1. Estimate of the error.} Similarly as in the proof of the Theorem \ref{thm:main}, we first note that since $(r, U)$ solves \eqref{sys:loc_euler_per}, it also solves
    \begin{equation}\label{sys:loc_euler_strong_per}
    \begin{aligned}
        \partial_tr + \DIV_x(r\, U) &= 0 &&\quad\mbox{in }(0,T)\times \T^d,\\
        \partial_t (r\,U)  + \DIV_x (r\, U \otimes U ) + r\nabla_x(L_{m,\varepsilon}^{per}\star(L_{m,\varepsilon}^{per}\star r)^{\gamma - 1}) + \kappa \,r \,U  &= e &&\quad\mbox{in }(0,T)\times \T^d, \\
    \end{aligned}
    \end{equation}
    for an error
    $$
        e := r\nabla_x(L_{m,\varepsilon}^{per}\star(L_{m,\varepsilon}^{per}\star r)^{\gamma - 1}) - r\nabla_x r^{\gamma - 1}.
    $$
We claim that
    \begin{align}\label{ineq:error_estimate_per}
       \left|\frac{e}{r}\right|\leq C(\gamma, L_m, \|r\|_{L^\infty(0, T; W^{2,\infty}(\T^d))}, r_*) \, \eps.
    \end{align}
Indeed, we have
    \begin{equation*}
    \begin{split}
         \left|\frac{e}{r}\right|&\leq |\nabla_x(L_{m,\varepsilon}^{per}\star(L_{m,\varepsilon}^{per}\star r)^{\gamma - 1}) - \nabla_x(L_{m,\varepsilon}^{per}\star(L_{m,\varepsilon}^{per}\star r^{\gamma - 1}))| + |\nabla_x(L_{m,\varepsilon}^{per}\star(L_{m,\varepsilon}^{per}\star r^{\gamma - 1})) - \nabla_x r^{\gamma-1}|\\
         &\leq (\gamma - 1)\int_{\T^d}L_{m,\eps}^{per}(y)|L_{m,\eps}^{per}\star\nabla_x r(x-y)||(L_{m,\eps}^{per}\star r)^{\gamma - 2}(x-y) - r^{\gamma- 2}(x-y)|\diff y\\
         &\quad +|\nabla_x(L_{m,\varepsilon}^{per}\star(L_{m,\varepsilon}^{per}\star r^{\gamma - 1})) - \nabla_x r^{\gamma-1}|.
    \end{split}
    \end{equation*}
Note that $r \geq r_*$ and so, $L_{m,\eps}^{per}\star r \geq r_*$.  The function $x \mapsto x^{\gamma - 2}$ is Lipschitz on $[r_*, \infty)$, hence using periodicity of $r$ we may estimate with a constant $C = C(\gamma, r_*)$
    \begin{align*}
    &|(L_{m,\eps}^{per}\star r)^{\gamma - 2}(x-y) - r^{\gamma- 2}(x-y)| \leq C\, |L_{m, \eps}^{per}\star r(x-y) - r(x-y)|\\
    &\leq C \int_{\T^d}\sum_{l\in\I^d}L_{m, \eps}(z+l)|r(x-y-z) - r(x-y)|\diff z  = C \int_{\T^d}\sum_{l\in\I^d}L_{m, \eps}(z+l)|r(x-y-z-l) - r(x-y)|\diff z \\
    &\leq C\, \|\nabla_x r\|_\infty  \int_{\T^d}\sum_{l\in\I^d}L_{m, \eps}(z+l)|z+l|\diff z  =C\, \eps \,  \|\nabla_x r\|_\infty \int_{\R^d}L_m(z)|z|\diff z.
    \end{align*}
    This implies
    \begin{align*}
(\gamma - 1)\int_{\T^d}L_{m,\eps}^{per}(y)|L_{m,\eps}^{per}\star\nabla_x r(x-y)||(L_{m,\eps}^{per}\star r)^{\gamma - 2}(x-y) - r^{\gamma- 2}(x-y)|\diff y \leq C\, \eps \, \|L_m(z)|z| \|_{1}\, \|\nabla_x r\|_\infty.
    \end{align*}
    For the other term with argue similarly to deduce
$$
        |\nabla_x(L_{m,\varepsilon}^{per}\star(L_{m,\varepsilon}^{per}\star r^{\gamma - 1})) - \nabla_x r^{\gamma-1}| \leq C\, \eps\,\|\nabla_x^2 r^{\gamma - 1}\|_\infty\leq C\eps(\|\nabla_x r\|_\infty^2 + \|r\|_\infty^{\gamma - 2}\|\nabla_x^2 r\|_\infty).
$$
This concludes the proof of \eqref{ineq:error_estimate_per}.\\

\noindent \underline{Step 2. The relative entropy inequality.} We claim that
    \begin{equation}\label{ineq:relative_entropy_3_per}
        \begin{split}
            \mathcal{E}^\varepsilon_{rel}(t) + \frac{1}{2}\mathrm{tr}(\mu^\varepsilon_t(\T^d))  - \mathcal{E}^\varepsilon_{rel}(0) &\leq -\kappa \int_0^t\int_{\T^d}|u^\varepsilon - U|^2\diff\rho^\varepsilon_\tau(x)\diff\tau + \int_0^t\int_{\T^d}(u^\varepsilon - U)\otimes (U - u^\varepsilon) : \nabla_x U\diff\rho^\varepsilon_\tau(x)\diff\tau\\
            &\quad + \int_0^t\int_{\T^d}U\cdot \nabla_x(L_{m,\eps}^{per}\star((L^{per}_{m,\eps}\star \diff\rho^\eps)^{\gamma - 1} - (L_{m,\varepsilon}^{per}\star r)^{\gamma - 1}))\diff\rho^\varepsilon_\tau(x)\diff\tau\\
            &\quad - \int_0^t\int_{\T^d}\nabla_x((\gamma - 1)L_{m,\eps}^{per}\star(\diff\rho^\eps - \diff r)(L_{m,\eps}^{per}\star r)^{\gamma - 2})L_{m,\eps}^{per}\star(r\,U)\, \diff x\diff\tau\\
            &\quad -\int_0^t\int_{\T^d}(u^\varepsilon - U)\cdot\frac{e}{r}\diff\rho^\varepsilon_\tau(x)\diff\tau - \int_0^t\int_{\T^d}\nabla_x U : \diff\mu^\varepsilon_\tau(x)\diff\tau.
        \end{split}
    \end{equation}
The derivation is analogous to the one given in the proof of the Theorem \ref{thm:main}. Nevertheless we give the full argument for the clarity of presentation. First, we test \eqref{eq:def_weak_cont_per} by $\frac{1}{2}|U|^2$ and $L_{m,\eps}^{per}\star(L_{m,\varepsilon}^{per}\star r)^{\gamma - 1}$ to obtain respectively
    \begin{equation}\label{eq:relative_U_2_per}
        \begin{split}
            \int_{\T^d}\frac{1}{2}|U|^2(t, x)\diff\rho^\varepsilon_t(x) - \int_{\T^d}\frac{1}{2}|U_0|^2\diff \rho^\varepsilon_0(x) = \int_0^t\int_{\T^d} \big(U\cdot \partial_t U +  u^\varepsilon\cdot(\nabla_x U)^{\top} U\big)\diff\rho^\varepsilon_\tau(x)\diff\tau
        \end{split}
    \end{equation}
    and
    \begin{equation}\label{eq:relative_K_r_per}
        \begin{split}
            &\int_{\T^d}(L_{m,\varepsilon}^{per}\star r)^{\gamma - 1}(t, x)(L_{m,\varepsilon}^{per}\star\diff\rho^\varepsilon)(t, x)\diff x - \int_{\T^d}(L_{m,\varepsilon}^{per}\star r_0)^{\gamma - 1}(L_{m,\varepsilon}^{per}\star\diff\rho^\varepsilon_0)\diff x\\
            &= \int_0^t\int_{\T^d} \big(\partial_t(L_{m,\eps}^{per}\star(L_{m,\varepsilon}^{per}\star r)^{\gamma - 1}) +  u^\varepsilon\cdot\nabla_x(L_{m,\eps}^{per}\star(L_{m,\varepsilon}^{per}\star r)^{\gamma - 1})\big)\diff\rho^\varepsilon_\tau(x)\diff\tau.
        \end{split}
    \end{equation}
    Furthermore, testing \eqref{eq:def_weak_moment_per} by $U$ gives us
    \begin{equation}\label{eq:relative_U_per}
        \begin{split}
            &\int_{\T^d} u^\varepsilon\cdot U\diff\rho^\varepsilon_t - \int_{\T^d}u^\varepsilon_0\cdot U_0\diff\rho^\varepsilon_0 = \int_0^t\int_{\T^d}\big(u^\varepsilon\cdot\partial_t U +  u^\varepsilon\otimes u^\varepsilon : \nabla_x U - \kappa \,u^\varepsilon\cdot U\big)\diff\rho_\tau^\varepsilon(x)\diff\tau\\
            &\qquad \qquad- \int_0^t\int_{\T^d}\nabla_x(L_{m,\eps}^{per}\star(L_{m,\varepsilon}^{per}\star \diff\rho^\varepsilon)^{\gamma - 1})\cdot U\diff \rho^\varepsilon_\tau(x)\diff\tau + \int_0^t\int_{\T^d}\nabla_x U:\diff\mu^\varepsilon_\tau(x)\diff \tau.
        \end{split}
    \end{equation}
    By adding together \eqref{eq:relative_U_2_per} and the energy inequality \eqref{ineq:def_weak_energy_per}, and subtracting from the result \eqref{eq:relative_U_per} as well as \eqref{eq:relative_K_r_per}, we get (note that in the first line we included an additional term $\frac{\gamma - 1}{\gamma}\int_{\T^d}(L^{per}_{m, \varepsilon}\star r)^\gamma\diff x - \frac{\gamma-1}{\gamma}\int_{\T^d}(L^{per}_{m, \varepsilon}\star r_0)^\gamma\diff x$ to form $\mathcal{E}^\varepsilon_{rel}(t)  - \mathcal{E}^\varepsilon_{rel}(0)$)
    \begin{equation}\label{ineq:relative_entropy_1_per}
        \begin{split}
            \mathcal{E}^\varepsilon_{rel}(t)+ \frac{1}{2}\mathrm{tr}(\mu^\varepsilon_t(\T^d))&  - \mathcal{E}^\varepsilon_{rel}(0) \leq -\kappa \int_0^t\int_{\T^d}|u^\varepsilon|^2\diff\rho^\varepsilon_\tau(x)\diff\tau + \frac{\gamma - 1}{\gamma}\int_{\T^d}(L^{per}_{m, \varepsilon}\star r)^\gamma\diff x - \frac{\gamma-1}{\gamma}\int_{\T^d}(L^{per}_{m, \varepsilon}\star r_0)^\gamma\diff x \\
            &\quad +\int_0^t\int_{\T^d} \big(U\cdot \partial_t U + u^\varepsilon\cdot(\nabla_x U)^{\top} U\big)\diff\rho^\varepsilon_\tau(x)\diff\tau\\
            &\quad - \int_0^t\int_{\T^d} \big(u^\varepsilon\cdot\partial_t U + u^\varepsilon\otimes u^\varepsilon : \nabla_x U - \kappa \, u^\varepsilon\cdot U - \nabla_x(L_{m,\eps}^{per}\star(L_{m,\varepsilon}^{per}\star \diff\rho^\eps)^{\gamma - 1})\cdot U\big)\diff \rho^\varepsilon_\tau(x)\diff\tau\\
            &\quad - \int_0^t\int_{\T^d}\nabla_x U : \diff\mu^\varepsilon_\tau(x)\diff\tau\\
            &\quad -\int_0^t\int_{\T^d}\big(\partial_t(L_{m,\eps}^{per}\star(L_{m,\varepsilon}^{per}\star r)^{\gamma - 1}) +  u^\varepsilon\cdot\nabla_x(L_{m,\eps}^{per}\star(L_{m,\varepsilon}^{per}\star r)^{\gamma - 1})\big)\diff\rho^\varepsilon_\tau(x)\diff\tau.
        \end{split}
    \end{equation}
    Now, since $r > 0$, from \eqref{sys:loc_euler_strong_per} we deduce that the strong velocity $U$ satisfies the equation
    $$
    \partial_t U + (U\cdot\nabla_x )U + \kappa\,U + \nabla_x (L_{m,\eps}^{per}\star(L_{m,\varepsilon}^{per}\star r)^{\gamma - 1}) = \frac{e}{r},
    $$
    thus we can write some terms appearing in \eqref{ineq:relative_entropy_1_per} as follows:
    \begin{equation*}
        \begin{split}
            &\int_0^t\int_{\T^d} \big(U\cdot \partial_t U +  u^\varepsilon\cdot(\nabla_x U)^{\top} U\big)\diff\rho^\varepsilon_\tau(x)\diff\tau - \int_0^t\int_{\T^d}\big( u^\varepsilon\cdot\partial_t U + u^\varepsilon\otimes u^\varepsilon : \nabla_x U\big)\diff\rho^\varepsilon_\tau(x)\diff \tau \\
            & = \int_0^t\int_{\T^d}\big((U - u^\varepsilon)\partial_t U + (U - u^\varepsilon)\cdot\nabla_x UU + (u^\varepsilon - U)\otimes (U - u^\varepsilon):\nabla_x U\big)\diff\rho^\varepsilon_\tau(x)\diff\tau\\
            & = \int_0^t\int_{\T^d} \Big( \kappa \, (u^\varepsilon - U)\cdot U + (u^\varepsilon - U)\cdot\nabla_x (L_{m,\eps}^{per}\star(L_{m,\varepsilon}^{per}\star r)^{\gamma - 1}) - (u^\varepsilon - U)\cdot\frac{e}{r}\Big)\diff \rho^\varepsilon_\tau(x)\diff\tau\\
            &\qquad + \int_0^t\int_{\T^d}(U - u^\varepsilon)\otimes (u^\varepsilon - U):\nabla_x U\diff\rho^\varepsilon_\tau(x)\diff\tau.
        \end{split}
    \end{equation*}
    Combining with \eqref{ineq:relative_entropy_1_per}, we obtain
    \begin{equation}\label{ineq:relative_entropy_2_per}
        \begin{split}
            \mathcal{E}^\varepsilon_{rel}(t)+ \frac{1}{2}\mathrm{tr}(\mu^\varepsilon_t(\T^d))&  - \mathcal{E}^\varepsilon_{rel}(0) \leq -\kappa \int_0^t\int_{\T^d}|u^\varepsilon - U|^2\diff\rho^\varepsilon_\tau(x)\diff\tau\\
            &\quad+ \frac{\gamma - 1}{\gamma}\int_{\T^d}(L^{per}_{m, \varepsilon}\star r)^\gamma\diff x - \frac{\gamma - 1}{\gamma}\int_{\T^d}(L^{per}_{m, \varepsilon}\star r_0)^\gamma\diff x \\
            &\quad+  \int_0^t\int_{\T^d}U\cdot\nabla_x (L_{m,\eps}^{per}\star(L_{m,\varepsilon}^{per}\star \diff\rho^\eps)^{\gamma - 1})\diff\rho^\varepsilon_\tau(x)\diff\tau\\
            &\quad +\int_0^t\int_{\T^d}\big((u^\varepsilon - U)\cdot\nabla_x (L_{m,\eps}^{per}\star(L_{m,\varepsilon}^{per}\star r)^{\gamma - 1}) + (u^\varepsilon - U)\otimes (U - u^\varepsilon) : \nabla_x U\big)\diff\rho^\varepsilon_\tau(x)\diff\tau\\
            &\quad -\int_0^t\int_{\T^d}(u^\varepsilon - U)\cdot\frac{e}{r}\diff\rho^\varepsilon_\tau(x)\diff\tau - \int_0^t\int_{\T^d}\nabla_x U : \diff\mu^\varepsilon_\tau(x)\diff\tau\\
            &\quad -\int_0^t\int_{\T^d}\big(\partial_t(L_{m,\eps}^{per}\star(L_{m,\varepsilon}^{per}\star r)^{\gamma - 1}) +  u^\varepsilon\cdot\nabla_x(L_{m,\eps}^{per}\star(L_{m,\varepsilon}^{per}\star r)^{\gamma - 1})\big)\diff\rho^\varepsilon_\tau(x)\diff\tau.
        \end{split}
    \end{equation}
    Let us focus on the terms with the kernel $L^{per}_{m, \varepsilon}$. Observe that
    \begin{align*}
        \frac{\gamma - 1}{\gamma}\int_{\T^d}(L^{per}_{m, \varepsilon}\star r)^\gamma\diff x - \frac{\gamma - 1}{\gamma}\int_{\T^d}(L^{per}_{m, \varepsilon}\star r_0)^\gamma\diff x = \int_0^t\int_{\T^d}\partial_t(L^{per}_{m,\varepsilon}\star r)^{\gamma - 1}L_{m,\eps}^{per}\star r\diff x\diff\tau.
    \end{align*}
    Thus, we can write some terms appearing in \eqref{ineq:relative_entropy_2_per} as 
    \begin{equation}\label{eq:relative_reduce_K_r_per}
        \begin{split}
            &\frac{\gamma - 1}{\gamma}\int_{\T^d}(L^{per}_{m, \varepsilon}\star r)^\gamma\diff x - \frac{\gamma - 1}{\gamma}\int_{\T^d}(L^{per}_{m, \varepsilon}\star r_0)^\gamma\diff x\\
            &\quad-\int_0^t\int_{\T^d}\big(\partial_t(L_{m,\eps}^{per}\star(L_{m,\varepsilon}^{per}\star r)^{\gamma - 1}) +  u^\varepsilon\cdot\nabla_x(L_{m,\eps}^{per}\star(L_{m,\varepsilon}^{per}\star r)^{\gamma - 1})\big)\diff\rho^\varepsilon_\tau(x)\diff\tau\\
            &\quad +\int_0^t\int_{\T^d}U\cdot\nabla_x (L_{m,\eps}^{per}\star(L_{m,\varepsilon}^{per}\star \diff\rho^\eps)^{\gamma - 1})\diff\rho^\varepsilon_\tau(x)\diff\tau\\
            &\quad +\int_0^t\int_{\T^d}(u^\varepsilon - U)\cdot\nabla_x (L_{m,\eps}^{per}\star(L_{m,\varepsilon}^{per}\star r)^{\gamma - 1})\diff\rho^\varepsilon_\tau(x)\diff\tau\\
            &= -\int_0^t\int_{\T^d} \partial_t(L_{m,\eps}^{per}\star(L_{m,\varepsilon}^{per}\star r)^{\gamma - 1})\diff(\rho^\varepsilon_\tau - r_\tau)(x)\diff\tau\\
            &\qquad+\int_0^t\int_{\T^d}U\cdot \nabla_x(L_{m,\eps}^{per}\star((L^{per}_{m,\eps}\star \diff\rho^\eps)^{\gamma - 1} - (L_{m,\varepsilon}^{per}\star r)^{\gamma - 1}))\diff\rho^\varepsilon_\tau(x)\diff\tau.
        \end{split}
    \end{equation}
    Using the continuity equation for the strong density $r$, and integrating by parts we have
    \begin{align*}
        &-\int_0^t\int_{\T^d} \partial_t(L_{m,\eps}^{per}\star(L_{m,\varepsilon}^{per}\star r)^{\gamma - 1})\diff(\rho^\varepsilon_\tau - r_\tau)(x)\diff\tau\\
        &= -\int_0^t\int_{\T^d}(\gamma - 1)L_{m,\eps}^{per}\star(\diff\rho^\eps - \diff r)(L_{m,\eps}^{per}\star r)^{\gamma - 2}L_{m,\eps}^{per}\star\p_t r\, \diff x\diff\tau\\
        &= \int_0^t\int_{\T^d}(\gamma - 1)L_{m,\eps}^{per}\star(\diff\rho^\eps - \diff r)(L_{m,\eps}^{per}\star r)^{\gamma - 2}L_{m,\eps}^{per}\star(\DIV_x(r \, U))\, \diff x\diff\tau\\
        &= -\int_0^t\int_{\T^d}\nabla_x((\gamma - 1)L_{m,\eps}^{per}\star(\diff\rho^\eps - \diff r)(L_{m,\eps}^{per}\star r)^{\gamma - 2})L_{m,\eps}^{per}\star(r\,U)\, \diff x\diff\tau.
    \end{align*}
    Combined with \eqref{eq:relative_reduce_K_r_per} this reduces \eqref{ineq:relative_entropy_2_per} to \eqref{ineq:relative_entropy_3_per}, as desired. \\

\noindent \underline{Step 3. Proof of \eqref{eq:rel_energy_before_Gronwall_general_pressure}.} 
We estimate the terms appearing in \eqref{ineq:relative_entropy_3_per}. The first term is included in \eqref{eq:rel_energy_before_Gronwall_general_pressure}. For the second term, is easily bounded since
    \begin{align}\label{ineq:relative_bound_1_per}
    \left|\int_0^t\int_{\T^d}(u^\varepsilon - U)\otimes (U - u^\varepsilon) : \nabla_x U\diff\rho^\varepsilon_\tau(x)\diff\tau\right|\leq \|\nabla_x U\|_\infty\int_0^t\int_{\T^d}|u^\varepsilon - U|^2\diff\rho^\varepsilon_\tau(x)\diff\tau \leq 2 \|\nabla_x U\|_\infty  \int_0^t \mathcal{E}^\varepsilon_{rel}(\tau) \diff \tau.
    \end{align}
    For the fifth term we apply the inequality on the error term \eqref{ineq:error_estimate_per}, the fact that $\rho^\eps$ is a probability measure and Young's inequality to see that
    \begin{equation}\label{ineq:relative_bound_4_per}
    \begin{split}
        \left|\int_0^t\int_{\T^d}(u^\varepsilon - U)\cdot\frac{e}{r}\diff\rho^\varepsilon_\tau(x)\diff\tau\right| &\leq \frac{1}{2}\int_0^t\int_{\T^d}|u^\varepsilon - U|^2\diff\rho^\varepsilon_\tau(x)\diff\tau + \frac{1}{2}\int_0^t\int_{\T^d}\left|\frac{e}{r}\right|^2\diff\rho^\varepsilon_\tau(x)\diff\tau\\
        &\leq \frac{1}{2} \int_0^t \mathcal{E}^\varepsilon_{rel}(\tau) \diff \tau + \eps^2 \,  C(\gamma, m, \|r\|_{L^\infty(0, T; W^{2,\infty}(\T^d))}, r_*).
    \end{split}
    \end{equation}
    Furthermore, for the term involving the concentration measure $\mu$ we can bound using Proposition \ref{prop:bound_symmetric_measure_trace}
    \begin{align}\label{ineq:relative_bound_5_per}
        \left|\int_0^t\int_{\T^d}\nabla_x U : \diff\mu^\varepsilon_\tau(x)\diff\tau\right| \leq C\|\nabla_x U\|_\infty\int_0^t\mathrm{tr}(\mu^\varepsilon_\tau(\T^d))\diff\tau.
    \end{align}
    Hence, we only need to discuss the two terms with kernel $L_{m,\varepsilon}$. We split them by introducing the commutator error to remove the nonlocal operator from $U$ as follows
    \begin{equation*}
        \begin{split}
            &\int_0^t\int_{\T^d}U\cdot \nabla_x(L_{m,\eps}^{per}\star((L^{per}_{m,\eps}\star \diff\rho^\eps)^{\gamma - 1} - (L_{m,\varepsilon}^{per}\star r)^{\gamma - 1}))\diff\rho^\varepsilon_\tau(x)\diff\tau\\
            &\quad - \int_0^t\int_{\T^d}\nabla_x((\gamma - 1)L_{m,\eps}^{per}\star(\diff\rho^\eps - \diff r)(L_{m,\eps}^{per}\star r)^{\gamma - 2})L_{m,\eps}^{per}\star(rU) \diff x\diff\tau\\
            &=\int_0^t\int_{\T^d}(L_{m,\eps}^{per}\star(U\diff\rho^\eps) - (L^{per}_{m,\eps}\star\diff\rho^\eps)\, U)\cdot\nabla_x((L^{per}_{m,\eps}\star \diff\rho^\eps)^{\gamma - 1} - (L_{m,\varepsilon}^{per}\star r)^{\gamma - 1})\diff x\diff\tau\\
            &\quad +\int_0^t\int_{\T^d} (L^{per}_{m,\eps}\star\diff\rho^\eps)\, U\cdot\nabla_x((L^{per}_{m,\eps}\star \diff\rho^\eps)^{\gamma - 1} - (L_{m,\varepsilon}^{per}\star r)^{\gamma - 1})\diff x\diff\tau\\
            &\quad -\int_0^t\int_{\T^d}(L_{m,\eps}^{per}\star(rU) - (L_{m, \eps}^{per}\star r)\, U)\cdot\nabla_x((\gamma - 1)L_{m,\eps}^{per}\star(\diff\rho^\eps - \diff r)(L_{m,\eps}^{per}\star r)^{\gamma - 2})\diff x\diff\tau\\
            &\quad -\int_0^t\int_{\T^d} (L_{m, \eps}^{per}\star r)\, U\cdot\nabla_x((\gamma - 1)L_{m,\eps}^{per}\star(\diff\rho^\eps - \diff r)(L_{m,\eps}^{per}\star r)^{\gamma - 2})\diff x\diff\tau =:I_1 + I_2 - I_3 - I_4.
        \end{split}
    \end{equation*}
We first handle $I_2-I_4$. We rewrite it by integrating by parts
    \begin{equation*}
        \begin{split}
            I_2 - I_4 &= -\int_0^t\int_{\T^d}\DIV_x U((L_{m,\eps}^{per}\star\diff\rho^\eps)^{\gamma} - (L_{m, \eps}^{per}\star\diff\rho^\eps)\,(L_{m,\eps}^{per}\star r)^{\gamma - 1})\diff x\diff\tau\\
            &\quad +(\gamma - 1)\int_0^t\int_{\T^d}\DIV_x U((L_{m,\eps}^{per}\star r)^{\gamma - 1}\, (L_{m, \eps}^{per}\star\diff\rho^\eps) - (L_{m,\eps}^{per}\star r)^{\gamma})\diff x\diff\tau\\
            &\quad -\int_0^t\int_{\T^d}U\cdot\left(\nabla_x\left(\frac{1}{\gamma}(L_{m,\eps}^{per}\star\diff\rho^\eps)^{\gamma}\right) - (L_{m,\eps}\star r)^{\gamma - 1}\nabla_xL_{m,\eps}^{per}\star\diff\rho^\eps\right)\diff x\diff\tau\\
            &\quad +\int_0^t\int_{\T^d}U\cdot\left(\nabla_x(L_{m,\eps}^{per}\star r)^{\gamma - 1}L_{m,\eps}^{per}\star\diff\rho^\eps - \frac{\gamma - 1}{\gamma}\nabla_x(L_{m,\eps}^{per}\star r)^{\gamma}\right)\diff x\diff\tau\\
            &=\left(\frac{1}{\gamma} - 1\right)\int_0^t\int_{\T^d}\DIV_x U\, d_{\gamma}(L_{m,\eps}^{per}\star\diff\rho^\eps,L_{m,\eps}^{per}\star r)\diff x\diff\tau.
        \end{split}
    \end{equation*}
    Hence,
    \begin{equation}\label{ineq:relative_bound_2_per}
        \begin{split}
            |I_2 - I_4|\leq C\|U\|_{L^\infty(0, T; W^{1, \infty}(\T^d))}\int_0^t\int_{\T^d}d_{\gamma}(L_{m,\eps}^{per}\star\diff\rho^\eps,L_{m,\eps}^{per}\star r)\diff x\diff\tau.
        \end{split}
    \end{equation}
    We move on to $I_1-I_3$. Recall that in Section~\ref{sect:auxilliary_results_periodic_case} we introduced
    \begin{align*}
        \mathcal{A}(x,y) = \nabla_x L_{m, \varepsilon}^{per}(x-y)\cdot (U(x) - U(y)) + \DIV_x U(x)\,  L_{m,\eps}^{per}(x-y).
    \end{align*}
Then, integrating by parts gives us
\begin{equation*}
\begin{split}
I_1-I_3
&= -\int_0^t\int_{\T^d}
\DIV_x\left(
L_{m,\eps}^{per}\star(U\diff\rho^\eps)
-(L_{m,\eps}^{per}\star\diff\rho^\eps)\,U
\right)
\left(
(L_{m,\eps}^{per}\star\diff\rho^\eps)^{\gamma-1}
-(L_{m,\eps}^{per}\star r)^{\gamma-1}
\right)\diff x\diff\tau\\
&\quad+\int_0^t\int_{\T^d}
\DIV_x\left(
L_{m,\eps}^{per}\star(rU)
-(L_{m,\eps}^{per}\star r)\,U
\right)
(\gamma-1)(L_{m,\eps}^{per}\star r)^{\gamma-2}
L_{m,\eps}^{per}\star(\diff\rho^\eps-\diff r)
\diff x\diff\tau\\
&=\int_0^t\int_{\T^d}\int_{\T^d}
\mathcal A(x,y)\diff\rho^\eps(y)
\left(
(L_{m,\eps}^{per}\star\diff\rho^\eps)^{\gamma-1}
-(L_{m,\eps}^{per}\star r)^{\gamma-1}
\right)(x) \diff x\diff\tau\\
&\quad-\int_0^t\int_{\T^d}\int_{\T^d}
\mathcal A(x,y)r(y)
(\gamma-1)(L_{m,\eps}^{per}\star r)^{\gamma-2}
L_{m,\eps}^{per}\star(\diff\rho^\eps-\diff r)(x)
\diff y\diff x\diff\tau.
\end{split}
\end{equation*}
 Then, recalling $d_{\gamma-1}(f,g) = f^{\gamma-1} - (\gamma-1)\, g^{\gamma-2}\, f + (\gamma-2)\, g^{\gamma-1}$, we observe that that factors containing $L_{m,\eps}^{per}$ above sum up to $d_{\gamma - 1}(L_{m, \varepsilon}^{per}\star\diff\rho^\eps, L_{m, \varepsilon}^{per}\star r)$. Therefore, writing $\diff\rho^\eps(y) = \diff\rho^\eps(y) - \diff r(y) + \diff r(y)$ in the first integral gives us
    \begin{equation}\label{eq:I_1_I_2_form}
    \begin{split}
        I_1 - I_3 &= \int_0^t\int_{\T^d}\int_{\T^d}\mathcal{A}(x,y)(\diff\rho^\eps(y) - \diff r(y))\, d_{\gamma - 1}(L_{m, \varepsilon}^{per}\star\diff\rho^\eps, L_{m, \varepsilon}^{per}\star r)(x)\diff x\diff\tau\\
        &\quad + \int_0^t\int_{\T^d}\int_{\T^d}\mathcal{A}(x,y)(\diff\rho^\eps(y) - \diff r(y))(\gamma - 1)(L_{m, \varepsilon}^{per}\star r)^{\gamma - 2}L_{m, \varepsilon}^{per}\star (\diff\rho^\eps - \diff r)(x)\diff x\diff\tau\\
        &\quad + \int_0^t\int_{\T^d}\int_{\T^d}\mathcal{A}(x,y)r(y)d_{\gamma - 1}(L_{m, \varepsilon}^{per}\star\diff\rho^\eps, L_{m, \varepsilon}^{per}\star r)(x)\diff y\diff x\diff\tau  =: H_1 + H_2 - H_3.
    \end{split}
\end{equation}
For the term $H_1$ we decompose using Lemma \ref{lem:operator-mathcal-A-decomposition}
\begin{equation*}
    \begin{split}
        H_1  & = \sum_{j = 1}^d\int_{0}^t\int_{\T^d}\sum_{i = 1}^d\partial_{x_j}U_i(x)\, P_{ij} \ast (\diff\rho^\eps - \diff r)  \,d_{\gamma - 1}(L_{m, \varepsilon}^{per}\star\diff\rho^\eps, L_{m, \varepsilon}^{per}\star r)(x)\diff x\diff\tau\\
        & \qquad+ \int_0^t\int_{\T^d}\int_{\T^d}S_2(x,y)\, (\diff\rho^\eps(y) - \diff r(y))\, d_{\gamma - 1}(L_{m, \varepsilon}^{per}\star\diff\rho^\eps, L_{m, \varepsilon}^{per}\star r)(x)\diff y\diff x\diff\tau  =: H_{1,1} + H_{1,2}.
    \end{split}
\end{equation*}
Using H\"{o}lder's inequality, \eqref{eq:bound_P_ij}, and Lemma \ref{lem:d_gamma_inequalities}
\begin{equation}\label{ineq:x_1_1_bound}
\begin{split}
    |H_{1,1}|&\leq \|U\|_{L^\infty(0, T;W^{1,\infty}(\T^d))}\int_0^t\sum_{i, j}\|P_{ij}\star(\diff\rho^\eps - \diff r)\|_{L^\gamma(\T^d)}\|d_{\gamma - 1}(L_{m, \varepsilon}^{per}\star\diff\rho^\eps, L_{m, \varepsilon}^{per}\star r)\|_{L^{\gamma/(\gamma - 1)}(\T^d)}\diff\tau\\
    &\leq C\|U\|_{L^\infty(0, T;W^{1,\infty}(\T^d))}\int_0^t\|L_{m,\eps}^{per}\star (\diff\rho^\eps - \diff r)\|_{L^\gamma(\T^d)}\|d_{\gamma - 1}(L_{m, \varepsilon}^{per}\star\diff\rho^\eps, L_{m, \varepsilon}^{per}\star r)\|_{L^{\gamma/(\gamma - 1)}(\T^d)}\diff\tau\\
    &\leq C\|U\|_{L^\infty(0, T;W^{1,\infty}(\T^d))}\int_0^t\int_{\T^d}d_\gamma(L_{m, \varepsilon}^{per}\star\diff\rho^\eps, L_{m, \varepsilon}^{per}\star r)(x)\diff x\diff\tau.
\end{split}
\end{equation}
Moreover, by \eqref{eq:bound_on_the_remainder_S2}, H\"{o}lder's inequality, Young's inequality, Lemma \ref{lem:d_gamma_inequalities} and the energy bounds \eqref{ineq:def_weak_energy_per} we can estimate
\begin{equation}\label{ineq:x_1_2_bound}
    \begin{split}
        |H_{1,2}|&\leq C\eps\|U\|_{L^\infty(0, T;W^{2,\infty}(\T^d))}\int_0^t\int_{\T^d}L_{m+4,\eps}^{per}\star(\diff\rho^\eps + \diff r)(x)\, d_{\gamma - 1}(L_{m, \varepsilon}^{per}\star\diff\rho^\eps, L_{m, \varepsilon}^{per}\star r)(x)\diff x\diff\tau\\
        &\leq C\eps\|U\|_{L^\infty(0, T;W^{2,\infty}(\T^d))}\int_0^t\|L_{4,\eps}^{per}\star L_{m,\eps}^{per}\star(\diff\rho^\eps + \diff r)\|_{L^\gamma(\T^d)}\|d_{\gamma - 1}(L_{m, \varepsilon}^{per}\star\diff\rho^\eps, L_{m, \varepsilon}^{per}\star r)\|_{L^{\gamma/(\gamma - 1)}(\T^d)}\diff\tau\\
        &\leq C\eps\|U\|_{L^\infty(0, T;W^{2,\infty}(\T^d))}\int_0^t\|L_{m,\eps}^{per}\star(\diff\rho^\eps + \diff r)\|_{L^\gamma(\T^d)}\|d^{(\gamma - 1)/\gamma}_{\gamma}(L_{m, \varepsilon}^{per}\star\diff\rho^\eps, L_{m, \varepsilon}^{per}\star r)\|_{L^{\gamma/(\gamma - 1)}(\T^d)}\diff\tau\\
        &\leq C\eps\|U\|_{L^\infty(0, T;W^{2,\infty}(\T^d))}\int_0^t\int_{\T^d}(L_{m, \varepsilon}^{per}\star\diff\rho^\eps)^\gamma+(L_{m, \varepsilon}^{per}\star r)^\gamma\diff x\diff\tau\\
        &\leq C\eps\|U\|_{L^\infty(0, T; W^{2,\infty}(\T^d))}.
    \end{split}
\end{equation}
For the term $H_2$, we also decompose the integrand using Lemma \ref{lem:operator-mathcal-A-decomposition}
\begin{equation*}
    \begin{split}
        H_2 &= \sum_{j = 1}^d\int_{0}^t\int_{\T^d}\sum_{i = 1}^d\partial_{x_j}U_i(x)\, P_{ij}\ast(\diff\rho^\eps(y) - \diff r(y)) \,(\gamma - 1)(L_{m, \varepsilon}^{per}\star r)^{\gamma - 2}L_{m, \varepsilon}^{per}\star (\diff\rho^\eps - \diff r)(x)\diff y\diff x\diff\tau\\
        & \qquad+ \int_0^t\int_{\T^d}\int_{\T^d}S_2(x,y) (\diff\rho^\eps(y) - \diff r(y))(\gamma - 1)(L_{m, \varepsilon}^{per}\star r)^{\gamma - 2}L_{m, \varepsilon}^{per}\star (\diff\rho^\eps - \diff r)(x)\diff y\diff x\diff\tau\\
        &=:H_{2,1} + H_{2,2}.
    \end{split}
\end{equation*}
Then, using H\"{o}lder's inequality, \eqref{eq:bound_P_ij}, and Lemma \ref{lem:d_gamma_inequalities} (this step also uses that $r\geq r_*$) we obtain
\begin{equation}\label{ineq:x_2_1_bound}
    \begin{split}
        |H_{2,1}|&\leq (\gamma - 1)\|r\|^{\gamma - 2}_{L^\infty((0, T)\times \T^d)}\|U\|_{L^\infty(0, T; W^{1,\infty}(\T^d))}\sum_{i,j}\int_0^t\|P_{ij}\star(\diff\rho^\eps - \diff r)\|_{L^2(\T^d)}\|L_{m, \varepsilon}^{per}\star (\diff\rho^\eps - \diff r)\|_{L^2(\T^d)}\diff\tau\\
        &\leq C(\gamma - 1)\|r\|^{\gamma - 2}_{L^\infty((0, T)\times \T^d)}\|U\|_{L^\infty(0, T; W^{1,\infty}(\T^d))}\int_0^t\|L_{m, \varepsilon}^{per}\star (\diff\rho^\eps - \diff r)\|^2_{L^2(\T^d)}\diff\tau\\
        &\leq C(\gamma - 1)\|r\|^{\gamma - 2}_{L^\infty((0, T)\times \T^d)}\|U\|_{L^\infty(0, T; W^{1,\infty}(\T^d))}\int_0^t\int_{\T^d}d_\gamma(L_{m, \varepsilon}^{per}\star\diff\rho^\eps, L_{m, \varepsilon}^{per}\star r)(x)\diff x\diff\tau.
    \end{split}
\end{equation}
For $H_{2,2}$ we argue similarly to $H_{1,2}$ to obtain
\begin{align}\label{ineq:x_2_2_bound}
    |H_{2,2}|\leq C\,\eps\,\|r\|^{\gamma - 2}_{L^\infty((0, T)\times \T^d)}\|U\|_{L^\infty(0, T; W^{2,\infty}(\T^d))},
\end{align}
where the constant $C$ depends on the energy of the initial conditions due to the bounds \eqref{ineq:def_weak_energy_per}. For the third integral we again write
\begin{equation*}
    \begin{split}
        H_3&=\sum_{j = 1}^d\int_{0}^t\int_{\T^d}\sum_{i = 1}^d\partial_{x_j}U_i(x)\, (P_{ij} \ast r) \,d_{\gamma - 1}(L_{m, \varepsilon}^{per}\star\diff\rho^\eps, L_{m, \varepsilon}^{per}\star r)(x)\diff y\diff x\diff\tau\\
        & \qquad+ \int_0^t\int_{\T^d}\int_{\T^d}S_2(x,y)\,r(y)\,d_{\gamma - 1}(L_{m, \varepsilon}^{per}\star\diff\rho^\eps, L_{m, \varepsilon}^{per}\star r)(x)\diff y\diff x\diff\tau =:H_{3,1} + H_{3,2}.
    \end{split}
\end{equation*}
For $H_{3,1}$ we use the explicit definition of $P_{ij}$ (see \eqref{eq:definition_pij_different_ij}--\eqref{eq:definition_pij_same_ij} in the proof of Lemma~\ref{lem:operator-mathcal-A-decomposition}) and H\"{o}lder's inequality to estimate
\begin{equation}\label{ineq:x_3_1_aux_bound_1}
    \begin{split}
        &|H_{3,1}|\leq\\
        & C\|U\|_{L^\infty_t W^{1,\infty}_x}\sum_{i,j}\int_0^t\|R_i^{per}R_j^{per}((L^{per}_{2,\eps}\star L_{m, \varepsilon}^{per})\star r - L_{m, \varepsilon}^{per}\star r)\|_{L^\gamma(\T^d)}\|d_{\gamma - 1}(L_{m, \varepsilon}^{per}\star\diff\rho^\eps, L_{m, \varepsilon}^{per}\star r)\|_{L^{\gamma/(\gamma - 1)}(\T^d)}\diff\tau.
    \end{split}
\end{equation}
Here, $R_i^{per}$ is the Riesz transform on $\T^d$. Since $R_i^{per}$ are bounded operators from $L^\gamma(\T^d)$ to $L^\gamma(\T^d)$ we can estimate (using as well Young's inequality and the periodicity of $r$)
\begin{equation*}
    \begin{split}
        &\|R_i^{per}R_j^{per}((L^{per}_{2,\eps}\star L_{m, \varepsilon}^{per})\star r - L_{m, \varepsilon}^{per}\star r)\|^\gamma_{L^\gamma(\T^d)}\\
        &\leq C\|(L^{per}_{2,\eps}\star L_{m, \varepsilon}^{per})\star r - L_{m, \varepsilon}^{per}\star r\|^\gamma_{L^\gamma(\T^d)}\\
        &= C\int_{\T^d}\left|\int_{\T^d}L_{2,\eps}^{per}(y)((L_{m,\eps}^{per}\star r)(x-y) - (L^{per}_{m,\eps}\star r)(x))\diff y\right|^\gamma\diff x\\
        &\leq C\int_{\T^d}\left(\int_{\T^d}\sum_{l\in\I^d}L_{2,\eps}(y+l)|(L_{m,\eps}^{per}\star r)(x-y) - (L^{per}_{m,\eps}\star r)(x)|\diff y\right)^\gamma\diff x\\
        & = C\int_{\T^d}\left(\int_{\T^d}\sum_{l\in\I^d}L_{2,\eps}(y+l)|(L_{m,\eps}^{per}\star r)(x-l-y) - (L^{per}_{m,\eps}\star r)(x)|\diff y\right)^\gamma\diff x\\
        &\leq C\|\nabla(L^{per}_{m, \eps}\star r)\|^\gamma_{L^\infty(\T^d)}\left(\sum_{l\in\I^d}\int_{\T^d}L_{2,\eps}(y+l)|y+l|\diff y\right)^\gamma\\
        &\leq C\eps^\gamma\|\nabla_x r\|^\gamma_{L^\infty(\T^d)}\left(\int_{\R^d}|L_2(y) y|\diff y\right)^\gamma.
    \end{split}
\end{equation*}
Hence, coming back to \eqref{ineq:x_3_1_aux_bound_1}, using Lemma \ref{lem:d_gamma_inequalities} and the energy bounds \eqref{ineq:def_weak_energy_per} we get
\begin{equation}\label{ineq:x_3_1_bound}
    \begin{split}
        |H_{3,1}|&\leq C\eps\|U\|_{L^\infty(0, T; W^{1,\infty}(\T^d))}\int_0^t\|d_{\gamma - 1}(L_{m, \varepsilon}^{per}\star\diff\rho^\eps, L_{m, \varepsilon}^{per}\star r)\|_{L^{\gamma/(\gamma - 1)}(\T^d)}\diff\tau\\
        &\leq C\eps\|U\|_{L^\infty(0, T; W^{1,\infty}(\T^d))}\int_0^t\left(\int_{\T^d}(L_{m, \varepsilon}^{per}\star\diff\rho^\eps)^\gamma + (L_{m, \varepsilon}^{per}\star r)^\gamma \diff x\right)^{(\gamma - 1)/\gamma}\diff\tau\\
        &\leq C\eps \|U\|_{L^\infty(0, T; W^{1,\infty}(\T^d))}.
    \end{split}
\end{equation}
Lastly, we argue similarly to $H_{1,2}$ and $H_{2,2}$ to estimate
\begin{equation}\label{ineq:x_3_2_bound}
    \begin{split}
        |H_{3,2}|\leq C\,\eps\, \|r\|_\infty,
    \end{split}
\end{equation}
for a constant $C> 0$ dependant on the energy bounds \eqref{ineq:def_weak_energy_per}. From \eqref{ineq:x_1_1_bound}, \eqref{ineq:x_1_2_bound}, \eqref{ineq:x_2_1_bound}, \eqref{ineq:x_2_2_bound}, \eqref{ineq:x_3_1_bound} and \eqref{ineq:x_3_2_bound} we deduce an inequality for \eqref{eq:I_1_I_2_form} of the form 
\begin{equation}\label{ineq:relative_bound_3_per}
    |I_1 - I_3| \leq C_A \int_0^t\int_{\T^d}d_\gamma(L_{m, \varepsilon}^{per}\star\diff\rho^\eps, L_{m, \varepsilon}^{per}\star r)(x)\diff x\diff\tau +C_B \, \eps,
\end{equation}
where $C_A$, $C_B$ depend on $d, m, \gamma, r_*, \|r\|_{L^\infty(0, T; W^{1,\infty}(\T^d))}, \|U\|_{L^\infty(0, T;W^{2,\infty}(\T^d))}$, and $C_B$ depend additionally on $\|u_0^\eps\|_{L^2(\T^d; \diff\rho_0^\eps)}$ through the exploited energy bounds \eqref{ineq:def_weak_energy_per}. Combining \eqref{ineq:relative_bound_1_per}, \eqref{ineq:relative_bound_4_per}, \eqref{ineq:relative_bound_5_per} as well as \eqref{ineq:relative_bound_2_per}, \eqref{ineq:relative_bound_3_per} together with \eqref{ineq:relative_entropy_3_per} we obtain \eqref{eq:rel_energy_before_Gronwall_general_pressure}.\\

\noindent \underline{Step 4. Proof of \eqref{eq:relative_entropy_after_gronwall_nonlinear_press_k_pos}, \eqref{eq:wasserstein_rate_nonlinear_press}, \eqref{eq:flat_rate_nonlinear_press}, \eqref{eq:negSS_rate_nonlinear_press}, and \eqref{eq:relative_entropy_after_gronwall_nonlinear_press_k_neg}.} When $\kappa \geq 0$, we can apply Gr\"{o}nwall's inequality to obtain \eqref{eq:relative_entropy_after_gronwall_nonlinear_press_k_pos} from \eqref{eq:rel_energy_before_Gronwall_general_pressure}. On the other hand, if $\kappa < 0$, we estimate it by $ -\kappa \int_0^t\int_{\T^d}|u^\varepsilon - U|^2\diff \rho^\varepsilon_\tau(x)\diff\tau \leq -2\,\kappa \int_0^t \mathcal{E}^\varepsilon_{rel}(\tau)\diff \tau$. Applying Gr\"{o}nwall's inequality yields \eqref{eq:relative_entropy_after_gronwall_nonlinear_press_k_neg}. Now, having \eqref{eq:relative_entropy_after_gronwall_nonlinear_press_k_pos} (or \eqref{eq:relative_entropy_after_gronwall_nonlinear_press_k_neg} when $\kappa <0$), we have a control on three terms
$$
\frac{1}{2}\int_{\T^d}|u^\varepsilon - U|^2\diff\rho^\varepsilon_t(x) + \frac{1}{\gamma}\int_{\T^d}d_\gamma(L^{per}_{m, \varepsilon}\star\diff\rho^\varepsilon, L_{m,\eps}^{per} \star\diff r)\diff x + \frac{1}{2}\mathrm{tr}(\mu^\varepsilon_t(\T^d)) =         \mathcal{E}_{rel}^\varepsilon(t)+ \frac{1}{2}\mathrm{tr}(\mu^\varepsilon_t(\T^d)).
$$ 
To prove \eqref{eq:wasserstein_rate_nonlinear_press} and  \eqref{eq:flat_rate_nonlinear_press} it is now sufficient to apply Lemma \ref{lem:estimate_wasserstein_flat} and \eqref{eq:relative_entropy_after_gronwall_nonlinear_press_k_pos} (or \eqref{eq:relative_entropy_after_gronwall_nonlinear_press_k_neg} if $\kappa < 0$). Regarding \eqref{eq:negSS_rate_nonlinear_press}, by Proposition~\ref{prop:bound_symmetric_measure_trace}, we have
    $
    \|\mu^\varepsilon_t\|_{TV}\leq C\, \mathrm{tr}(\mu_t^\varepsilon(\T^d)).
    $
    To prove the estimate in $H^{-m}(\T^d)$ norm, we first note that by Lemma \ref{lem:d_gamma_inequalities},
$$
C(r_*, \|r\|_{\infty})\, \int_{\T^d} |L^{per}_{m, \varepsilon}\star\diff\rho^\varepsilon-  L_{m,\eps}^{per} \star\diff r|^2 \diff x \leq  \frac{1}{\gamma}\int_{\T^d}d_\gamma(L^{per}_{m, \varepsilon}\star\diff\rho^\varepsilon, L_{m,\eps}^{per} \star\diff r)\diff x \leq \mathcal{E}_{rel}^\varepsilon(t).
$$
The (LHS) controls $\|\rho^\varepsilon - r\|_{H^{-m}(\T^d)}$ (see, e.g., \eqref{eq:nonlocal_operator_gives_convergence_in_neg_SS_linear_pressure}) which concludes the proof of \eqref{eq:negSS_rate_nonlinear_press}. Finally, the convergence statements follow now directly from \eqref{eq:relative_entropy_after_gronwall_nonlinear_press_k_pos} (or \eqref{eq:relative_entropy_after_gronwall_nonlinear_press_k_neg} if $\kappa < 0$), \eqref{eq:wasserstein_rate_nonlinear_press}, \eqref{eq:flat_rate_nonlinear_press}, and \eqref{eq:negSS_rate_nonlinear_press}.

\section{Proof of Theorem \ref{thm:main_particles}}\label{sect:proof_main_result_for_particles}

\noindent We will need the following result whose variants are well-known in the literature, see, e.g., \cite[Theorem~1.41]{MR4309603}.

\begin{lem}\label{lem:particle_approx_of_a_gen_measure}
Let $r \in \mathcal{P}(\R^d)$ be such that $\int_{\R^d} |x|^{2+\delta} \diff r(x) < \infty$ for some $\delta>0$. Then, there exist explicitly given disjoint measurable sets $\{A_i\}_{i=1}^{k^d+1}$, with $\cup_{i=1}^{k^d+1} A_i = \R^d$, and points $\{x_i\}_{i=1}^{k^d+1}$ with $x_i \in A_i$ such that, defining
$
r^N := \sum_{i=1}^{k^d+1} r(A_i)\, \delta_{x_i},\, N = k^d+1,
$
we have
$$
W^2_2(r, r^N) \leq C\, k^{-\frac{2\,\delta}{(2+\delta)}},\quad N = k^d+1,
$$
where $C$ depends only on the dimension $d$ and $\int_{\R^d} |x|^{2+\delta} \diff r$, and $W_2$ denotes the 2-Wasserstein distance; see Appendix~\ref{sect:wasserstein_flat}. If $\R^d$ is replaced with $\T^d$, then the assumption $\int_{\R^d} |x|^{2+\delta} \diff r(x) < \infty$ is not needed and we have $W_2^2(r^N, r) \leq C\, k^{-2}$, $N = k^d$.
\end{lem}

\begin{rem}
One can obtain rate $C\, k^{-2}$ in the case of $\R^d$ at the cost that the points $x_i$ are not explicit by using techniques from the quantization theory, see \cite[Corollary 6.7]{MR1764176}. However, as our main interest in the current work is the numerics, we prefer explicit construction over the abstract result.
\end{rem}
\begin{proof}[Proof of Lemma \ref{lem:particle_approx_of_a_gen_measure}]
We first discuss the case of $\R^d$. We split $\R^d$ for the open cube $(-R,R)^d$ and $\R^d \setminus (-R,R)^d$. The cube $(-R,R)^d$ can be splitted into $k^d$ disjoint cubes $\{Q_i\}_{i=1}^{k^d}$ of side length $2R/k$ (not necessarily open or closed). We define $Q_{k^d+1} :=\R^d \setminus (-R,R)^d$ and we let $r^N := \sum_{i=1}^{k^d+1} r(Q_i)\, \delta_{x_i}$ with $x_i \in Q_i$ for $i = 1, ..., k^d$ and some $x_{k^d+1} \in Q_{k^d+1}$ satisfying $|x_{k^d+1}|=R$. We can estimate
\begin{align*}
W^2_2(r^N, r) &\leq \sum_{i=1}^{k^d+1}   \int_{Q_i} |x-x_i|^2 \diff r(x) \leq \sum_{i=1}^{k^d}\frac{4\,R^2\,d}{k^2} \, r(Q_i) + 4\, \int_{Q_{k^d+1}} |x|^2 \diff r(x) \\
&\leq \frac{4\,R^2\,d}{k^2} + 4 \, R^{-\delta} \, \int_{\R^d} |x|^{2+\delta} \diff r(x) \leq C \, (R^2 \, k^{-2} + R^{-\delta}),
\end{align*}
where on the last set we used that $|x-x_{k^d+1}|^2 \leq 2\,(|x|^2 + |x_{k^d+1}|^2) \leq 4\, |x|^2$. We choose $R$ such that $R^2 \, k^{-2} =  R^{-\delta}$, i.e. $R = k^{\frac{2}{2+\delta}}$ so that we obtain $W_2^2(r^N, r) \leq C\,  k^{-\frac{2\,\delta}{(2+\delta)}}$ for $N = k^d + 1$. On $\T^d$, we can just set $R=1$ and there is no need to consider the set $|x|\geq R$ so we obtain $W_2^2(r^N, r) \leq C\, k^{-2}$ with $N = k^d$. 
\end{proof}

\noindent We are now ready to prove Theorem \ref{thm:main_particles}.

\begin{proof}[Proof of Theorem \ref{thm:main_particles}]
We first discuss the case of $\R^d$. We want to apply Theorem \ref{thm:main} to the solution $\rho^{N,\varepsilon}$ (which is a solution by Proposition \ref{prop:particle_sln_is_a_sol}). To use the bounds from Theorem \ref{thm:main}, we need to estimate the initial kinetic energy $\sum_{i=1}^{k^d+1} |U^0(X_i(0)|^2 \, r(A_i)$ (it corresponds to the term $\|u_0^{\eps}\|_{L^2(\R^d;\diff\rho_0^{\eps})}$ in the formulation of Theorem~\ref{thm:main}) and $\mathcal{E}_{\mathrm{rel}}^\varepsilon(0) =  \int_{\R^d}|L_{m, \varepsilon}\star(\diff r^N_0 - \diff r_0)|^2\diff x$. For the initial kinetic energy, we use that $U^0$ is bounded and $\sum_{i=1}^{k^d+1} r(A_i) = 1$.  For the $\mathcal{E}_{\mathrm{rel}}^\varepsilon(0)$ term, by transforming to Fourier variables, we have
$$
\mathcal{E}_{\mathrm{rel}}^{N,\varepsilon}(0) = \int_{\R^d}  \frac{1}{(1 + 4\pi^2|\eps\,\xi|^2)^{m}} \big|\widehat{r_0} - \widehat{r_0^N}\big|^2 \diff \xi. 
$$ 
A direct computation shows that 
$$
|\widehat{r_0}(\xi) - \widehat{r_0^N}(\xi)| = \left|\int_{\R^d} e^{-2\pi i x \xi} \diff(r_0 - r_0^N)(x)\right| \leq 2\,\pi\, (1+|\xi|) \, W_1(r_0, r_0^N) \leq 2\,\pi\, (1+|\xi|) \, W_2(r_0, r_0^N),
$$ 
where we used that $W_1(r_0, r_0^N) \leq W_2(r_0, r_0^N)$, cf. \eqref{eq:monotonicity_Wasserstein}, and the dual representation \eqref{eq:W_1_dual_representation}. It follows that
$$
|\mathcal{E}_{\mathrm{rel}}^{N,\varepsilon}(0)| \leq 4\, \pi^2\, W^2_2(r_0, r_0^N)\int_{\R^d}  \frac{(1+|\xi|)^2}{(1 + 4\pi^2|\xi|^2)^{m}}  \diff \xi \leq \frac{8\, \pi^2\, W^2_2(r_0, r_0^N)}{\varepsilon^2}  \int_{\R^d}  \frac{(1+|\eps\,\xi|)^2}{(1 + 4\pi^2|\eps\, \xi|^2)^{m}}  \diff \xi \leq C\, \frac{8\, \pi^2\, W^2_2(r_0, r_0^N)}{\varepsilon^{2+d}},
$$
where $C$ is a numerical constant which exists because $2(m-1)>(d-1)$. From Lemma \ref{lem:particle_approx_of_a_gen_measure} we know that $
W^2_2(r_0, r^N_0) \leq C\, k^{-\frac{2\,\delta}{(2+\delta)}},\, N = k^d + 1,
$ so that we deduce 
$
|\mathcal{E}_{\mathrm{rel}}^{N,\varepsilon}(0)| \leq C\, {\varepsilon^{-(2+d)}}\, k^{-\frac{2\,\delta}{(2+\delta)}},\, N = k^d + 1, 
$
where $C$ depends only on the dimension $d$ and $\int_{\R^d} |x|^{2+\delta} \diff r$. Plugging these estimates to \eqref{eq:final_estimate_W_2_on_R_d}--\eqref{eq:final_estimate_neg_SS_on_R_d}, we obtain \eqref{eq:rate_SPH_Rd_W2}--\eqref{eq:rate_SPH_Rd_negSS}. \\

\noindent On $\mathbb{T}^d$, we want to apply Theorem \ref{thm:main_per}. As above, we need to estimate the initial kinetic energy and the term $\mathcal{E}_{\mathrm{rel}}^{N,\varepsilon}(0)$. The initial kinetic energy is controlled in the same way but the term $\mathcal{E}_{\mathrm{rel}}^{N,\varepsilon}(0)$ will require more work. We observe that
$
  \mathcal{E}_{rel}^{N,\varepsilon}(0) :=  \frac{1}{\gamma}\int_{\T^d}d_\gamma(L^{per}_{m, \varepsilon}\star r_0^N, L_{m,\eps}^{per} \star r_0)\diff x.
$
The function $d_{\gamma}$ satisfies
$$
d_\gamma(f, g) = f^\gamma - \gamma g^{\gamma - 1}f + (\gamma - 1)g^\gamma= f^{\gamma} - g^{\gamma} + \gamma \, g^{\gamma-1}\, (f-g) \leq C(\gamma)\, (f^{\gamma-1} + g^{\gamma-1}) \, |f-g|,
$$
where $C$ is a numerical constant depending on $\gamma$, as the function $x \mapsto x^{\gamma}$ is Lipschitz on $[0,M]$ with constant $\gamma\, M^{\gamma-1}$. Clearly,
$$
|L_{m,\eps}^{per} \star r_0^N|^{\gamma-1} \leq \|L_{m,\eps}^{per}\|_{L^{\infty}(\T^d)}^{\gamma-1} \, \|r_0^N\|^{\gamma-1}_{\mathcal{M}( \T^d)} \leq C(m,d)\, \eps^{-d},
$$
$$
|L_{m,\eps}^{per} \star r_0|^{\gamma-1} \leq \|L_{m,\eps}^{per}\|_{L^{1}(\T^d)}^{\gamma-1} \, \|r_0\|^{\gamma-1}_{L^{\infty}(\T^d)} \leq C(\gamma,   \|r_0\|_{L^{\infty}(\T^d)}),
$$
where $ \|r_0^N\|_{\mathcal{M}(\T^d)} = 1$ denotes the total mass of the measure $r_0^N$. It follows that 
$$
 \mathcal{E}_{rel}^{N,\varepsilon}(0) \leq C(\gamma, m, d,  \|r_0\|_{L^{\infty}(\T^d)}) \, \eps^{-d}\, \| L_{m,\eps}^{per} \star (r_0 - r_0^N)\|_{L^1(\T^d)}.
$$
Using \eqref{eq:W_1_dual_representation}, we estimate further
$$
\| L_{m,\eps}^{per} \star (r_0 - r_0^N)\|_{L^1(\T^d)} \leq |\T^d| \, \|\nabla L_{m,\eps}^{per}\|_{L^{\infty}(\T^d)} \,  W_2(r_0, r_0^N) \leq C(m,d)\, \eps^{-(d+1)} \, W_2(r_0, r_0^N).
$$ 
In order to estimate $\|\nabla L_{m,\eps}^{per}\|_{L^{\infty}(\T^d)}$ we need the assumption $m\geq d+1$. Indeed, by \eqref{bessel_function_at_0} and \eqref{bessel_potential_derivative}, if $m \geq d+1$ the function $\nabla L_{m}$ is bounded. Then, we have to control the series
$$
\nabla L_{m,\eps}^{per}(x) = \frac{1}{\eps^{d+1}} \sum_{l\in\I^d} (\nabla L_{m}) \left(\frac{x +l}{\eps} \right). 
$$
As $\eps$ is getting smaller, $\frac{x +l}{\eps}$ is larger in magnitude and the series can be uniformly controlled by using the exponential tail \eqref{bessel_function_at_infinity}. We conclude that 
$$
\mathcal{E}_{rel}^{N,\varepsilon}(0) \leq C(\gamma, m, d,  \|r_0\|_{L^{\infty}(\T^d)}) \,\eps^{-(2d+1)} \, W_2(r_0, r_0^N).
$$
From Lemma \ref{lem:particle_approx_of_a_gen_measure} we know that $
W_2(r_0, r^N_0) \leq C\, k^{-1}$, $N = k^d$. Plugging these estimates into \eqref{eq:wasserstein_rate_nonlinear_press}--\eqref{eq:negSS_rate_nonlinear_press}, we obtain \eqref{eq:rate_SPH_Td_W2}--\eqref{eq:rate_SPH_Td_negSS}.
\end{proof}

\section*{Acknowledgements}
\noindent Piotr Gwiazda and Jakub Skrzeczkowski was supported by National Science Center, Poland through project no. 2023/51/B/ST1/01546. Jakub Woźnicki was supported by National Science Center, Poland through project no. 2021/43/O/ST1/03031. The authors would like to thank José A. Carrillo for bringing the works of Di Lisio, Grenier, and Pulvirenti to their attention. This work was partially supported by the Simons Foundation grant (award no. SFI-MPS-T-Institutes-00010825) and from State Treasury funds as part of a task commissioned by the Minister of Science and Higher Education under the project “Organization of the Simons Semesters at the Banach Center - New Energies in 2026-2028” (agreement no. MNiSW/2025/DAP/491).”
\appendix

\section{Properties of the Bessel potential}

\noindent In this section, we recall from the book of Aronszajn and Smith \cite{aronszajn1961theory} some of the properties of the Bessel potential. Given a function $L_m$ such that
$$
\hat{L}_m(\xi) = \frac{1}{(1 + 4\pi^2|\xi|^2)^{m/2}}, \quad m > 0,
$$
it is given by a semi-explicit formula
$$
L_m(x) = \frac{1}{2^{\frac{d + m - 2}{2}}\pi^{d/2}\Gamma\left(\frac{m}{2}\right)}F_{\frac{d - m}{2}}(|x|)|x|^{\frac{m - d}{2}},
$$
where $F$ is a modified Bessel function of the third kind (cf \cite[Eq. (4.1)]{aronszajn1961theory}). Following \cite[Section 4]{aronszajn1961theory} we note that for any $m > 0$ the function $L_m$ is decreasing, everywhere positive and integrable. Moreover, using the convention
$$
f\sim g, \text { as }z\to a,\text { whenever }\lim_{z\to a}\frac{f(z)}{g(z)} = 1,
$$
as $|x| \to 0$
\begin{equation}\label{bessel_function_at_0} 
\begin{split}
    &F_\nu(|x|) \sim 2^{|\nu| -1}\Gamma(|\nu|)|x|^{-|\nu|},\,\nu\neq 0,\\
    &F_0(|x|)\sim \ln\left(\frac{1}{|x|}\right).
\end{split}
\end{equation}
and
\begin{equation}\label{bessel_potential_at_0}
\begin{split}
    &L_m(x) \sim \frac{\Gamma\left(\frac{d - m}{2}\right)}{2^m \pi^{d/2}\Gamma\left(\frac{m}{2}\right)}|x|^{m - d}, \quad m < d,\\
    &L_d(x) \sim \frac{1}{2^{d-1}\pi^{d/2}\Gamma\left(\frac{d}{2}\right)}\ln\left(\frac{1}{|x|}\right),\\
    &L_m(x) \sim \frac{\Gamma\left(\frac{d - m}{2}\right)}{2^m \pi^{d/2}\Gamma\left(\frac{m}{2}\right)}, \quad m > d,
\end{split}
\end{equation}
as well as when $|x|\to +\infty$
\begin{align}\label{bessel_function_at_infinity}
    F_\nu(|x|)\sim \left(\frac{\pi}{2|x|}\right)^{1/2}\exp(-|x|),
\end{align}
and
\begin{align}\label{bessel_potential_at_infinity}
L_m(x) \sim \frac{1}{2^{\frac{d+m -1}{2}}\pi^{\frac{d-1}{2}}\Gamma\left(\frac{m}{2}\right)}|x|^{\frac{m - d - 1}{2}}\exp(-|x|).
\end{align}
Moving further, denote by $|x| = r$, then, by \cite[Eq. (4.4)]{aronszajn1961theory}
\begin{align}\label{bessel_potential_derivative}
\frac{d}{dr}L_m(r) = -\frac{1}{2^{\frac{m+d - 2}{2}}\pi^{d/2}\Gamma\left(\frac{m}{2}\right)}r^{\frac{m-d}{2}}F_{\frac{d - m + 2}{2}}(r).
\end{align}
Having all of the above, we may now show a crucial lemma.
\begin{lem}\label{lem:bound_for_nabla_bessel_potential}
    Suppose $L_m$ is a Bessel potential with $m > d$. Then, there exists a constant $C > 0$, such that
    $$
    \frac{|\nabla_x L_m(x)||x|^2}{L_4\star L_m(x)} \leq C.
    $$
\end{lem}
\begin{proof}
    Note that the function
    $$
    f(x) :=  \frac{|\nabla_x L_m(x)||x|^2}{L_4\star L_m(x)}
    $$
    is continuous on $\R^d\setminus\{0\}$. Therefore, the only possibility of being unbounded is at limits in $0$ and $+\infty$. We will show that the limits 
    $$
     \lim_{|x|\to 0}\frac{|\nabla_x L_m(x)||x|^2}{L_4\star L_m(x)}, \quad\lim_{|x|\to +\infty} \frac{|\nabla_x L_m(x)||x|^2}{L_4\star L_m(x)},
    $$
    are bounded above, which will end the proof. Using \eqref{bessel_potential_derivative}
    $$
    \frac{|\nabla_x L_m(x)||x|^2}{L_4\star L_m(x)} \leq \frac{1}{L_4\star L_m(x)}\frac{1}{2^{\frac{m+d - 2}{2}}\pi^{d/2}\Gamma\left(\frac{m}{2}\right)}|x|^{\frac{m-d + 4}{2}}F_{\frac{d - m + 2}{2}}(|x|) =: g(x).
    $$
    Then, by \eqref{bessel_potential_at_0} whenever $|x|\to 0$
    $$
    \frac{1}{L_4\star L_m(x)} = \frac{1}{L_{m+4}(x)}\sim \frac{2^{m + 4} \pi^{d/2}\Gamma\left(\frac{m+4}{2}\right)}{\Gamma\left(\frac{d - m-4}{2}\right)}
    $$
    and by \eqref{bessel_function_at_0}
    $$
    |x|^{\frac{m-d + 4}{2}}F_{\frac{d - m + 2}{2}}(|x|) \sim \left\{{\setstretch{2.0}\begin{array}{lll}
         2^{\frac{d - m}{2}}\Gamma\left(\frac{d - m + 2}{2}\right)|x|^{m - d + 1}, &&0 < m < d+2,\\
         \frac{|x|}{\ln\left(\frac{1}{|x|}\right)},&& m = d+2,\\
         2^{\frac{m - d - 4}{2}}\Gamma\left(\frac{m - d - 2}{2}\right)|x|^{3},&& m > d+2.
    \end{array}}\right.
    $$
    Hence,
    \begin{align*}
        \lim_{|x|\to 0}g(x) = 0.
    \end{align*}
    Moreover, by \eqref{bessel_potential_at_infinity} whenever $|x|\to +\infty$
    $$
    \frac{1}{L_4\star L_m(x)} = \frac{1}{L_{m+4}(x)}\sim 2^{\frac{d+m +3}{2}}\pi^{\frac{d-1}{2}}\Gamma\left(\frac{m+4}{2}\right)|x|^{\frac{d - 3 - m}{2}}\exp(|x|)
    $$
    and by \eqref{bessel_function_at_infinity}
    $$
    |x|^{\frac{m-d + 4}{2}}F_{\frac{d - m + 2}{2}}(|x|) \sim \left(\frac{\pi}{2}\right)^{1/2}|x|^{\frac{m - d + 3}{2}}\exp(-|x|).
    $$
    Thus,
    $$
    \lim_{|x|\to +\infty}g(x) = A,
    $$
    for some constant $A > 0$. This ends the proof.
\end{proof}

\section{Relative entropy estimates}

\begin{lem}\label{lem:d_gamma_inequalities}
    Suppose $\gamma \geq 2$, and define
    $$
    d_\gamma(x, y) = x^\gamma - \gamma y^{\gamma - 1}x + (\gamma - 1) y^\gamma, \quad x,y\in \R^+.
    $$
    Then, there exists a constant $1 > c > 0$ such that for any $x, y\in \R^+$
    \begin{enumerate}
        \item $$
        |x - y|^\gamma \leq d_\gamma(x, y),
        $$
        \item $$
        c\, d_{\gamma - 1}(x, y)^{\frac{\gamma}{\gamma  - 1}} \leq d_\gamma(x, y),
        $$
    \end{enumerate}
    Moreover, if $0 < y_0 \leq y\leq Y_1 < +\infty$, then
    \begin{enumerate}[start = 3]
        \item $$
        c|x - y|^2 \leq d_\gamma(x, y).
        $$
    \end{enumerate}
\end{lem}
\begin{proof}
    We begin with the first inequality. Notice that whenever $y = 0$, the inequality is trivial, hence, we may assume $y > 0$ and divide both sides by $y^\gamma$, and consider instead the inequality
    $$
    0\leq z^\gamma - \gamma z + (\gamma - 1) - |z - 1|^\gamma =: f(z),
    $$
    for $z = x/y$ and a function $f$ defined on $(0, +\infty)$. Note that $f(1) = 0$. We will show that for $0\leq z < 1$ we have $f'(z)\leq 0$ and for $z > 1$ we get $f'(z) \geq 0$, which will end the proof. For $z\in [0, 1)$
    $$
    f'(z) = \gamma(z^{\gamma - 1} - 1 + (1 - z)^{\gamma - 1}).
    $$
    Since $1 - z > 0$, we have
    $$
    1 = (1- z + z)^{\gamma - 1}\geq z^{\gamma - 1} + (1-z)^{\gamma - 1},
    $$
    which is what we needed. Similarly one shows the assertion for $z > 1$. For the second inequality, we first note that for $\gamma = 2$ it is trivial, and then for $\gamma > 2$ we similarly as before reduce the problem to the inequality
    $$
    c \leq \frac{z^\gamma - \gamma z  + (\gamma - 1)}{(z^{\gamma - 1} + (\gamma - 1)z + (\gamma - 2))^{\gamma/(\gamma - 1)}} =: g(z),
    $$
    for $z = x/y$ and $g$ defined on $(0, +\infty)$. Note that the numerator of $g(z)$ is equal to $0$ only at $z = 1$, hence $g$ has a possibility of reaching $0$ only at $z = 1$ or at the limit $z\to +\infty$. It is easy to check that
    $$
    \lim_{z\to +\infty}g(z) = 1.
    $$
    Meanwhile, at $z \to 1$, we use Taylor's expansion with Lagrange's remainder to deduce that there exist $\zeta_z, \xi_z\in [z, 1]$ such that
    $$
    \lim_{z\to 1}g(z) = \lim_{z\to 1}\frac{\gamma(\gamma - 1)\zeta_z^{\gamma - 2}(z-1)^2}{((\gamma - 1)(\gamma - 2)\xi_z^{\gamma - 3}(z-1)^2)^{\gamma/(\gamma - 1)}} = +\infty,
    $$
    which ends the proof. For the third inequality we refer to \cite[Remark 3.5]{lattanzio2017fromgas}.
\end{proof}

\section{Wasserstein and flat metrics}\label{sect:wasserstein_flat}

\noindent We briefly recall the definitions of the two metrics on spaces of measures that will be used throughout the paper: the $p$-Wasserstein distance (with $p\in [1,\infty)$) and the flat metric. Let $\Omega=\R^d$ or $\T^d$. For two probability measures $\mu_1,\mu_2$ with finite $p$-th moments (i.e. $\int_{\Omega} |x|^p \diff \mu_1, \int_{\Omega} |x|^p \diff \mu_2 < \infty$), the $p$-Wasserstein distance is defined by 
\begin{equation}\label{eq:def_Wasserstein_p}
W_p^p(\mu_1,\mu_2)
=
\min_{\gamma\in\Pi(\mu_1,\mu_2)}
\int_{\Omega\times\Omega}|x-y|^p\diff\gamma(x,y),
\end{equation}
where $\Pi(\mu_1,\mu_2)$ denotes the set of couplings of $\mu_1$ and $\mu_2$, namely
$$
\Pi(\mu_1,\mu_2)
=
\left\{
\gamma\in\mathcal P(\Omega\times\Omega):
(\pi_1)_\#\gamma=\mu_1,\;
(\pi_2)_\#\gamma=\mu_2
\right\},
$$
where $\pi_1$ and $\pi_2$ denote the canonical projections onto the first and second components, respectively. We will also use the flat metric, which is defined for finite Radon measures and, in particular, does not require the two measures to have the same mass. For $\mu,\nu\in\mathcal M(\Omega)$, we set
\begin{equation}\label{eq:def_flat_metric}
d_f(\mu,\nu)
=
\sup\left\{
\int_\Omega \varphi\diff(\mu-\nu):
\varphi\in C^1(\Omega),\;
\|\varphi\|_\infty\leq1,\;
\|\nabla_x\varphi\|_\infty\leq1
\right\}.
\end{equation}
It is known (see, e.g., \cite[Theorem 1.82]{MR4309603}) that $W_1$ admits a similar characterization
\begin{equation}\label{eq:W_1_dual_representation}
W_1(\mu,\nu) 
=
\sup\left\{
\int_\Omega \varphi\diff(\mu-\nu):
\varphi\in C^1(\Omega),\;
\|\nabla_x\varphi\|_\infty\leq1
\right\}.
\end{equation}
Moreover, by applying Holder inequality for the definition \eqref{eq:def_Wasserstein_p}, we obtain
\begin{equation}\label{eq:monotonicity_Wasserstein}
W_p(\mu, \nu) \leq W_q(\mu, \nu) \qquad 1 \leq p \leq q.
\end{equation}
We refer to \cite{MR3409718} and \cite{MR4309603} for a comprehensive treatment of Wasserstein and flat metrics, respectively.

\begin{lem}\label{lem:estimate_wasserstein_flat} Suppose that $(r,U)$ and $(\rho,u)$ satisfy the continuity equations 
$$
\partial_t r + \DIV_x(r\, U)=0, \qquad \partial_t \rho_t + \DIV_x (\rho\, u)=0,
$$
with $r_t,\rho_t\in\mathcal P_2(\mathbb{R}^d)$ for all $t\in[0,T]$, $\nabla U\in L^1(0,T;L^\infty(\mathbb{R}^d))$, $\int_0^T \int_{\R^d} |U_t|^2 \diff r_t(x) \diff t, \int_0^T \int_{\R^d} |u_t|^2 \diff \rho_t(x) \diff t < \infty$. Then, for every $t\in[0,T]$, 
\begin{equation}\label{eq:general_formula_wasserstein_between_slns_cont_eq}
W_2^2(\rho_t,r_t) \leq \exp\left( \int_0^t \bigl(1+2\|\nabla U_s\|_{L^\infty(\R^d)}\bigr)\diff s \right) \left( W_2^2(\rho_0,r_0) + \int_0^t\int_{\mathbb{R}^d} |u_s-U_s|^2\diff \rho_s \diff s \right).
\end{equation}
Additionally, if $U \in L^{\infty}(0,T; W^{1,\infty}(\R^d))$, then
\begin{equation}\label{eq:general_formula_flat_between_momentum_cont_eq}
d_f(\rho_t u_t,r_t U_t)
\leq
\left(\int_{\R^d}|u_t-U_t|^2\diff\rho_t\right)^{\frac{1}{2}}
+ C_U
\left(
W_2^2(\rho_0,r_0)
+
\int_0^t\int_{\R^d}|u_s-U_s|^2\diff\rho_s\diff s
\right)^{1/2},
\end{equation}
where $C_U
:=
\|U\|_{L^\infty(0,T;W^{1,\infty}(\R^d))}
\exp\left(
\int_0^T
\left(\frac{1}{2}+\|\nabla U_s\|_{L^\infty(\R^d)}\right)\diff s
\right)$. The same is true when $\R^d$ is replaced with $\T^d$.
\end{lem} 
\begin{proof} Let $\pi_t$ be an optimal coupling in the Kantorovich problem between $\rho_t$ and $r_t$. By the differentiation formula for the squared Wasserstein distance along absolutely continuous curves in $\mathcal P_2(\mathbb{R}^d)$ (see \cite[Theorem~8.4.7]{AGS}), for a.e. $t\in(0,T)$, 
$$
\frac{1}{2}
\frac{\diff}{\diff t}W_2^2(\rho_t,r_t) = \int_{\mathbb{R}^d\times\mathbb{R}^d} (x-y)\cdot\bigl(u_t(x)-U_t(y)\bigr)\diff \pi_t(x,y).
$$
Writing $ u_t(x)-U_t(y) = u_t(x)-U_t(x)+U_t(x)-U_t(y) $ and using the Lipschitz continuity of $U_t$, we obtain 
\begin{align*} \frac{1}{2} \frac{\diff}{\diff t}W_2^2(\rho_t,r_t) &\leq \int_{\mathbb{R}^d\times\mathbb{R}^d} |x-y|\,|u_t(x)-U_t(x)|\diff \pi_t(x,y)  +\|\nabla U_t\|_{L^\infty(\R^d)} \, W_2^2(\rho_t,r_t) \\ &\leq W_2(\rho_t,r_t) \left( \int_{\mathbb{R}^d}|u_t-U_t|^2\diff \rho_t \right)^{1/2} +\|\nabla U_t\|_{L^\infty(\R^d)} \, W_2^2(\rho_t,r_t). \end{align*}
By Cauchy-Schwarz inequality, 
$$
\frac{\diff}{\diff t} W_2^2(\rho_t,r_t) \leq \bigl(1+2\,\|\nabla U_t\|_{L^\infty(\R^d)}\bigr)\, W_2^2(\rho_t,r_t) + \int_{\mathbb{R}^d}|u_t-U_t|^2\diff \rho_t.
$$
Hence, \eqref{eq:general_formula_wasserstein_between_slns_cont_eq} follows from Gronwall's inequality. To prove \eqref{eq:general_formula_flat_between_momentum_cont_eq}, we estimate
    \begin{align*}
        d_f(\rho\, u, r\,U) \leq d_f(\rho\, u, \rho\, U) + d_f(\rho\,U, r\,U) \leq \left(\int_{\R^d}|u - U|^2\diff\rho(x)\right)^{\frac{1}{2}} + \|U\|_{L^\infty(0, T; W^{1,\infty}(\R^d))}\,d_f(\rho, r).
    \end{align*}
    In order to conclude the proof, it remains to notice that $d_f(\rho, r) \leq W_1(\rho, r)$ (which is a direct consequence of comparing \eqref{eq:def_Wasserstein_p} and \eqref{eq:def_flat_metric}) and $W_1(\rho,r) \leq W_2(\rho,r)$ cf. \eqref{eq:monotonicity_Wasserstein}.
    \end{proof}

\bibliographystyle{abbrv}
\bibliography{parpde_mo_discmodul}
\end{document}